\documentclass[leqno,12pt]{article}
\usepackage{amsmath, amsthm, mathrsfs}
\usepackage{amsfonts}
\usepackage{amssymb}
\usepackage{graphicx}
\usepackage{color}
\usepackage{hyperref}
\usepackage{ifpdf} 
\usepackage{geometry}
\newtheorem{theorem}{Theorem}[section]
\newtheorem{corollary}[theorem]{Corollary}
\newtheorem{lemma}[theorem]{Lemma}
\newtheorem{proposition}{Proposition}[section]

\newcommand{\re}{\mathbb{R}}
\newcommand{\ren}{{\mathbb{R}^N}}
\newcommand{\ve}{\varepsilon}

\newcommand{\var}{\varphi}
\newcommand{\la}{\lambda}

\newcommand{\dx}{\,{\rm d}x}

\newcommand{\dt}{\,{\rm d}t}

\numberwithin{equation}{section}
\def\qed{\,\unskip\kern 6pt \penalty 500
\raise -2pt\hbox{\vrule \vbox to8pt{\hrule width 6pt
\vfill\hrule}\vrule}\par}
\newcommand{\bc}{\color{blue}}

\definecolor{ao}{rgb}{0.0, 0.0, 1.0}

\newcommand{\nc}{\normalcolor}

\definecolor{darkblue}{rgb}{0.05, .05, .65}
\definecolor{darkgreen}{rgb}{0.1, .65, .1}
\definecolor{darkred}{rgb}{0.8,0,0}

\title{\bf Growing Solutions of the fractional\\ $p$-Laplacian equation in the \\Fast Diffusion Range
}
\author{ \sc Juan Luis V\'azquez\\ Univ. Aut\'onoma de Madrid }

\date{ 28 Feb 2021}

\begin{document}
\maketitle

\begin{abstract}

We establish existence, uniqueness as well as  quantitative estimates for  solutions
$u(t,x)$ to the fractional nonlinear diffusion equation, $\partial_t u
+{\mathcal L}_{s,p} (u)=0$, where ${\mathcal L}_{s,p}=(-\Delta)_p^s$ is the standard fractional $p$-Laplacian operator. We work in the  range of exponents $0<s<1$ and $1<p<2$, and in some sections $sp<1$. The equation is posed in the whole space $x\in {\mathbb R}^N$. We first obtain  weighted global integral estimates that allow establishing the existence of solutions for a class of large data that is proved to be roughly optimal. We use the estimates to study the class of self-similar solutions of forward type, that we describe in detail when they exist. We also explain what happens when possible self-similar solutions do not exist.  We establish the dichotomy positivity versus extinction for nonnegative solutions at any given time. We analyze the conditions for extinction in  finite time.
\end{abstract}


\section{Introduction}\label{se.intro}

Nonlocal integro-differential operators and equations are  attracting increasing attention because of their mathematical interest, and also because of their multiple applications to fields like  Statistical Mechanics, Finance, Ecology, Image processing, Fluid Mechanics, and others. In particular, the fractional Laplacian operator is defined for suitable functions $u(x)$, $x\in {\mathbb R}^N$, as
	\begin{equation}\label{eq.frlap}
	(-\Delta)^s u(x) = c_{n,s}  PV \int_{{\mathbb R}^N} \frac{u(x) - u(y)}{|x-y|^{n+2s}} \,dy,
	\end{equation}
	where $s\in(0,1)$, $c_{n,s}>0$ is a normalization constant, and $PV$ denotes that the integral is taken in the principal value sense, cf. the classical references \cite{Riesz1949, Landkof66, Stein1970}, as well as recent ones like \cite{BucurBook16, CaffarelliSilv2007, GarofaloMR3916700, Kwasnicki2017, Ros-MR3447732}. Together with its variants, the operator plays a prominent role in describing anomalous diffusion, cf. \cite{VazTNLD2016} where local and nonlocal, linear and nonlinear models are compared.

We are concerned here with the nonlinear version given by the fractional $p$-Laplace operator
${\mathcal L}_{s,p}$ defined by the formula
\begin{equation}\label{frplap.op}
\displaystyle \qquad {\mathcal L}_{s,p}(u):= P.V.\int_{{\mathbb R}^N}\frac{\Phi_p(u(x,t)-u(y,t))}{|x-y|^{N+sp}}\,dy\,,
\end{equation}
where $1<p<\infty$, $\Phi_p(z)=|z|^{p-2}z,$  and $P.V.$ means principal value. To be precise, it can  be called the $s$-fractional $p$-Laplacian operator. It is  well-known from general theory that ${\mathcal L}_{s,p}$ is a maximal monotone operator in $L^2({\mathbb R}^N)$ with dense domain, more precisely the subdifferential of a convex Gagliardo functional.
Following  our previous papers \cite{Vazquez2020, VazFPL2-2020}, we will continue the study of the corresponding gradient flow, i.e., the evolution equation
\begin{equation}\label{frplap.eq}
\partial_t u + {\mathcal L}_{s,p} u=0,
\end{equation}
posed in the Euclidean space $x\in {\mathbb R}^N$, $ N\ge 1$, for $t>0$. We refer to it as the {\sl fractional $p$-Laplacian evolution equation}, FPLE for short.  Motivation and related equations for this model can be seen in the  \cite{Vazquez2020} and its references. There, the superlinear case $p>2$ was studied. The case $1<p< 2$, usually called fast diffusion, has been treated in \cite{VazFPL2-2020}, mainly in the class of solutions that are $L^1$-integrable with respect to the space variable. We pursue in this paper the analysis of  such a fractional fast diffusion for more general data, possibly growing at infinity, and devote much attention to settle different qualitative and quantitative issues, like positivity, self-similarity, mass conservation, and extinction.

\medskip

\noindent {\bf Outline of results.} As a starting novelty, the present paper introduces a basic tool called the \sl weighted $L^1$-estimate \rm that will play an important role in the existence theory for general classes of data. See whole details in Section \ref{sec.wl1est}.  Indeed, obtaining some weighted a priori $L^1$-estimate was the  motivation of the paper. The original idea goes back to the local nonlinear estimate proved by Herrero and Pierre in 1984 for the Fast Diffusion range of the Porous Medium Equation, cf. \cite{HerrPierre1984MR797051}. In their result the estimate had a stronger form, it has a local form  instead our weighted form. It was used by the authors to establish existence of solutions for the Cauchy problem in the whole space without any growth requirement on the initial data.

It took time to adapt that tool to fractional equations. Our estimate is based on the differential inequality
 \begin{equation}
|\frac{d}{dt}\int_{{\mathbb R}^N} u(t) \varphi \,dx| \le \big(\int_{{\mathbb R}^N} u(t) \varphi \,dx\big)^{p-1}\,K(\varphi),
\end{equation}
which is valid for all nonnegative solutions of the FPLE on the condition that $sp<1$ and $\var$ is a smooth positive function that decays in a controlled way as $|x|\to\infty$, roughly like a power between $|x|^{-(N+sp)}$ and $|x|^{-N-(sp/(p-1))}$.
Integration in time of this differential inequality gives the result  \eqref{est.wL1}, the version that is used in practice. It is also valid for differences of solutions $u=u_1-u_2\ge 0$. The proof is based on a kind of duality that has not been used before to our knowledge.
Our version is not as strong as the original one, and indeed it must be weaker since the derived existence results cannot be so general, as we will show below with precise counterexamples. Another close motivation for the present analysis comes from the study of the fractional version of the PME done in~\cite{BonfVaz2014MR3122168}.

Note that for $s\le 1/2$ we prove the estimate for the whole range of $p$ in fast diffusion,  $1<p<2$.  However, in the range $1/2<s<1$ we need the extra condition, $sp<1$, that seems to play a major role in the proof. The range is made clear in Figure~1.

\begin{figure}[t!]
 \centering
 \vspace{-2cm}
 \includegraphics[scale=0.25]{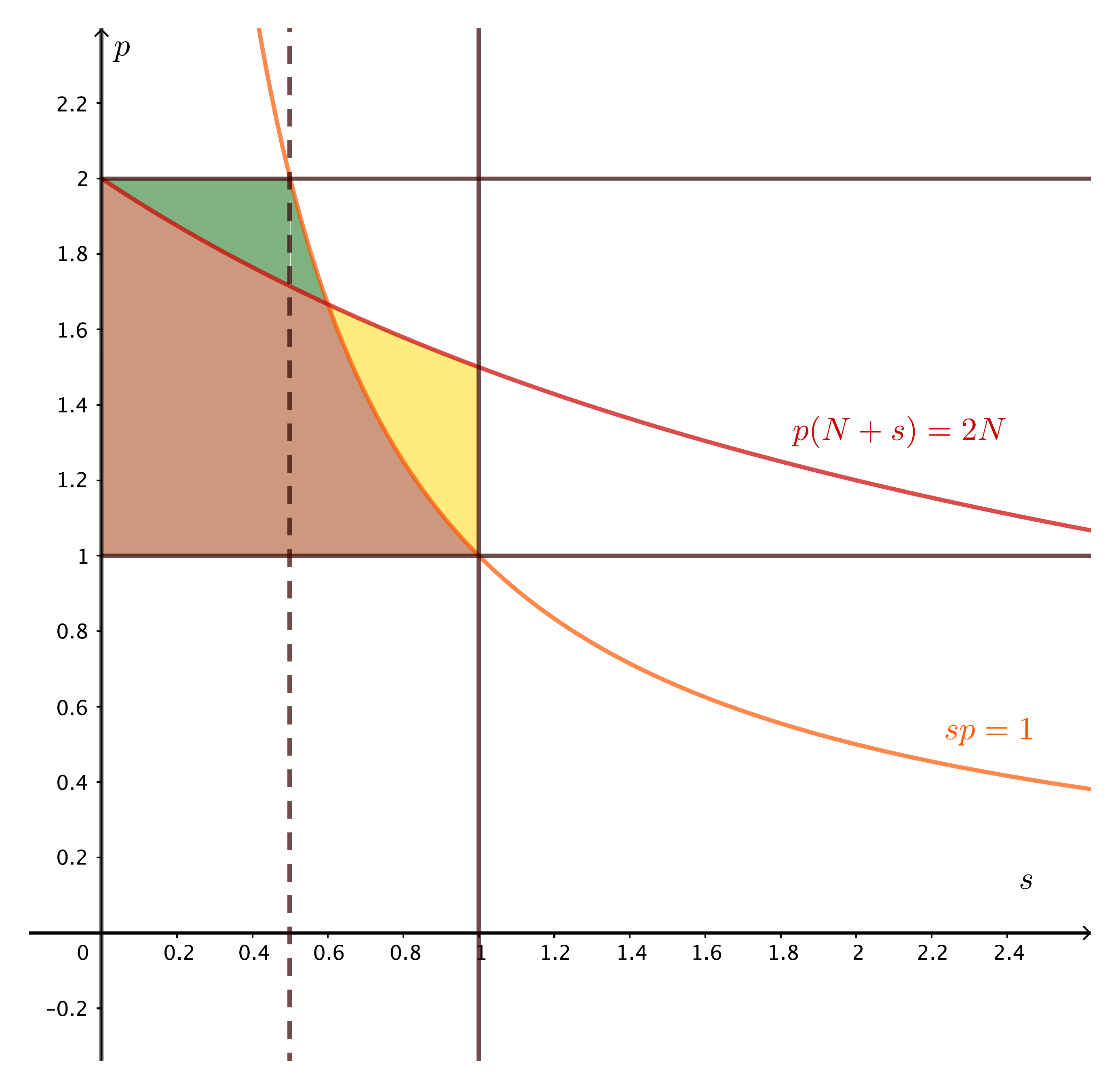}
 \caption{{Green and brown are the accepted regions.} White and green form the so-called Good Fast Diffusion range. {Yellow and brown correspond to Very Fast Diffusion.}}
  \label{fig:1}
\end{figure}

\medskip

\noindent $\bullet$
Section \ref{sec.ex} studies  the construction of solutions under a suitable growth limitation on the initial data. We prove existence of a finite ``candidate solution'' under almost optimal growth conditions, that we state here for brevity as
 \begin{equation}
u_0(x)=O(|x|^{\gamma}), \quad \mbox{with \ } \ \gamma<\frac{sp}{p-1}.
\end{equation}
The sharper condition is given in Formula \eqref{est.excond}. The solution is obtained by approximation from below with standard semigroup solutions (the ones studied in paper \cite{VazFPL2-2020} with initial data in $L^q$ spaces). The weighted $L^1$ estimate is then used to ensure that the limit of the approximate process produces  a finite locally integrable function $u(x,t)$. We prove that the limit is also a weak solution (in the sense of Definition \ref{def.ws}).  The constructed limit solutions are proved to be unique, independent of the chosen approximation. We use the name of  {\sl minimal weak solution} for the resulting function (often shortened to minimal solution). The precise results are stated in Theorem \ref{thm.exist} and \ref{thm.uniqlarge}. Solutions with changing sign are briefly mentioned.

The optimality of the power growth $\gamma< \gamma_1=sp/(p-1)$ in our existence condition is shown via the non-existence result  for the critical data $u_0(x)=|x|^{\gamma_1}$ that we prove in Theorem \ref{thm.instbu}. Note that in the limit $p=2$ we recover the known existence condition $\gamma<2s $ for the fractional heat equation, that was established in \cite{MR3211862, BonSirVazMR3614666}.

\medskip

\noindent $\bullet$ Section \ref{sec.masscon} is devoted to prove the remarkable dichotomy in behaviour for nonnegative solutions, we have \sl either positivity or extinction\rm. Precisely stated, for any nonnegative solution $u$ and any fixed time $t>0$, either $u(x,t)$ is positive for all $x$, or $u(x,t)=0$ for all $x\in{\mathbb R}^N$. In the last case, there exists a first time $T(u_0)>0$ where the solution vanishes identically, it is called the \sl extinction time\rm. Then, the solution must be continued for $t\ge T$ as  identically zero for all $x\in{\mathbb R}^N$. There is also the possibility $T(u_0)=\infty$, i.e., no extinction in finite time. In both  cases,  $u(x,t)$ is a locally strictly positive function of $x$  for any fixed $0<t<T$.

\medskip

\noindent $\bullet$ It is a fact that only a handful of explicit or semi-explicit solutions of the FPLE are known to date. Examples of known solutions are linear functions, fundamental solutions and very singular solutions, the two last classes being described in \cite{VazFPL2-2020}. This serves as motivation for the next part of the program, where we study the existence and properties of a large class of self-similar solutions.  In Section \ref{sec.self} we first review the scaling properties of the equation and analyze the possible self-similarities, that fall into three known types. Self-similar solutions of the standard type, called sometimes  forward self-similarity, take the form
 \begin{equation}
u(x, t)= t^{\gamma\beta} F(xt^{-\beta}).
\end{equation}
When these solutions exist they do for all time $t>0$ and preserve their size up to scaling in time.
They are  possible and they expand in space when the algebraic condition $sp+\gamma(2-p)>0$ is met, see for instance \cite{Vazquez2006}. This leads to the value of  $\beta=(sp+\gamma(2-p))^{-1}$ being positive.   In particular, we exclude two bizarre phenomena  for nontrivial data: there is no possible extinction or blow-up at a positive, finite time for this class.

We devote the next two sections to explore the existence of such self-similar solutions, always under the necessary condition $sp+\gamma(2-p)>0$ for forward expanding self-similarity. In Section \ref{sec.exss} we explore the existence of self-similar solutions with growing initial data of the form
\begin{equation}\label{init.power}
u_0(x)=A\,|x|^{\gamma},  \quad A,\gamma>0.
\end{equation}
We find the necessary limitation in the range of exponents for existence, $0< \gamma<\gamma_1=sp/(p-1)$, in agreement with the growth conditions of Section \ref{sec.ex}. We provide a  detailed description of the qualitative and quantitative properties of the family of solutions, see Theorem \ref{thm.ss.growing}. We have $u_t\ge 0$ for all $t>0,$ $x\in{\mathbb R}^N$. As we mentioned above, forward self-similarity implies no extinction in finite time.

A dramatic change occurs at the limit exponent $\gamma=\gamma_1$. Indeed, we prove an important result called {\sl instantaneous  blow-up}, cf.  Theorem \ref{thm.instbu}. Thus, if we try to obtain a limit solution by approximation from below with standard semigroup solutions,   as proposed in Section \ref{sec.ex}, the limit is infinite for all $t>0$ and $x\in \ren$. The same happens for all $\gamma>\gamma_1$ as an easy consequence of the maximum principle. Note that this means that $0<\gamma<\gamma_1=sp/(p-1)$ is a necessary and sufficient condition for the existence of self-similar solutions with growing data of power type $|x|^\gamma$.

We address in  Section \ref{sec.exss2} the question of existence of self-similar solutions for classes of decaying initial data, which take  the form \eqref{init.power} with $\gamma<0$. Solutions are taken in the limiting sense of Section \ref{sec.ex}. The theory offers no special novelties as long as $\gamma>-N$ so that $u_0$ is locally integrable, cf. Theorem \ref{decay.int}. A main difference is that now $u_t\le 0$ for all $t>0,$ $x\in\ren$.

There is a remaining interesting range where $sp+\gamma(2-p)>0 $ (so we can have forward self-similarity), and $\gamma\le -N$ (so that the initial function has a non-integrable singularity). This means that we are in the region where
$$
p_c<p<2, \qquad \mbox{and \ } \ -\frac{sp}{2-p}<\gamma\le -N.
$$
Here, $p_c=2N/(N+s)$ is the critical exponent that plays a big role in \cite{VazFPL2-2020}. This case is really curious. We find that there are two options:

(i) For $p_c<p<p_1$ we find a {\sl singular self-similar solution}, cf. Theorem \ref{ex.sing.sss}. The singularity stays at the origin for all times.

(ii) for $p_1\le p<2$ we find {\sl instantaneous blow-up}, hence no self-similar solution is possible, cf. Theorem \ref{nonex.sing.sss}.

Here, $p_1$ is the curious intermediate critical exponent identified in our previous paper \cite{VazFPL2-2020}, see Formula (1.10) and Figure 1 there. In \cite{VazFPL2-2020}, it appeared mainly as a separation of different asymptotic behaviour for the fundamental solutions. We have found here that it also plays an important role in the existence theory.

\medskip

\noindent $\bullet$ Next, we turn our attention to another appealing feature, {\sl extinction in finite time}. This is a classical topic of fast diffusion equations, see \cite{Vazquez2006}. The weighted a priori estimate  is very useful in studying extinction problems. We show that it can happen for $1<p<p_c$. We are able to establish fine conditions for existence and nonexistence of finite time extinction in Section \ref{sec.vfd}. The Morrey space $M^{q*}(\ren)$, with $q_*=N(2-p)/sp$, plays an important role, together with the Lebesgue space $L^{q*}(\ren)$. Constructed counter-examples show the sharpness of the results.

\bc The borderline case $p=p_c$ is very special since it separates two open ranges with very different qualitative behaviour. We contribute to the analysis of this case by proving the conservation of mass is still true as in the range $p>p_c$, see Theorem \ref{thm.mc.fple}. We provide a very nontrivial proof that uses a number of new tools. Previously known information is given in Section 16 of paper \cite{VazFPL2-2020}.\nc

 Naturally, ideas and techniques coming from the study of related non-fractional equations  are very useful in order to investigate the fractional model and to evaluate the results. In our case, we have to look at the standard  $p$-Laplacian equation,  $u_t=\mbox{div} (|\nabla u|^{p-2}\nabla u)$, i.e.,  the non-fractional case that corresponds to the limit $s=1$.  This model has been studied for fast diffusion, $1<p<2$, in a number of references,   \cite{BonVazFD2010, BonIagVaz2010, DiBenedetto93, DiBenHerrMR1066761, Lindqvist06, RazSanVaz2008, Vazquez2006}.   We will make comments in the paper on  the existing knowledge about the  standard equation  in relation to our results.

We devote one section to examine the connection of the paper to the work in the fractional Porous Medium Equation, in particular to paper  \cite{BonfVaz2014MR3122168} done in collaboration with M. Bonforte.
We conclude with an Appendix and a section on comments and open problems. \nc

\medskip

\noindent  {\bf Notations.} In the whole paper the assume the values $0<s<1$ and $1<p<2$; we will not recall it again as a general rule. The particular subrange of $p$ will be always carefully announced, in particular which sections depend on the assumption $sp<1$. We use the notation ${\mathcal L}_{s}\varphi=(-\Delta)^{s}\varphi$, omitting the exponent $p=2$ in ${\mathcal L}_{s,2}$.
We sometimes write $u(t)$ instead of $u(x,t)$ for convenience, when we want to stress that it is a function of $x$ with parameter $t$. In the sequel we often use the notation $a^m$ to mean the signed power $|a|^{m-1}a$ of a quantity $a\in \re$, $a\ne 0$.  All of this is done for brevity when no confusion is to be feared. We call good fast diffusion range the set of parameters where $p_c<p<2$, i.e., when $(N+s)p>2N$, and very fast diffusion range when $1<p<p_c$, we recall that $p_c=2N/(N+s)$. This distinction is quite important in the $L^1$ theory, in questions like mass conservation, existence of Barenblatt solutions, large-time behaviour and others, as shown in \cite{VazFPL2-2020}. It is relevant here in some results.

\section{The weighted $L^1$ estimate}\label{sec.wl1est}

We take the announced values of the parameters $s, p$  with the condition $sp<1$. We will prove an a priori estimate for all nonnegative solutions $u$ of the FPLE \eqref{frplap.eq} defined in a strip $Q=\ren\times (0,T)$.

\subsection{Presentation} In this section we work with semigroup solutions defined in all $L^q$ spaces and in later sections we extend the scope of the result  to the class of weak solutions by approximation. In order to state the main estimate we introduce that concept of weighted mass at time $0\le t<T$
\begin{equation}
X(t; u,\varphi)=\int_\ren u(t) \varphi \,dx\,,
\end{equation}
where the weight $\var$ is a positive function to be specified next.  First, we introduce the operator ${\mathcal M}_{s'}$  by the formula
\begin{equation}\label{frplap.opM}
({\mathcal M}_{s'}\var)(x):= P.V.\int_{\ren}\frac{|\var(x,t)-\var(y,t)|}{|x-y|^{N+2s'}}\,dy\,.
\end{equation}
We remark that when $0<2s'<1$ this operator is well-defined and bounded for bounded and uniformly Lipschitz continuous functions since the singularity at $x=y$ is integrable. Of course, we have ${\mathcal M}_{s'}(\var)\ge 0$ and also \ $|{\mathcal L}_{s'}( \var)|\le {\mathcal M}_{s'}(\var)$.

\medskip

\noindent{\bf The class $\mathcal C=\mathcal C(s,p)$.} The class of suitable weight functions for our main estimate is formed by the smooth and positive functions $\var$ defined in $\ren$ such that ${\mathcal M}_{sp/2}(\varphi)$ is locally bounded and
\begin{equation}\label{cond.phi}
C(\var):=\int \frac{|{\mathcal M}_{sp/2}\varphi(x)|^{1/(2-p)}}{\varphi(x)^{(p-1)/(2-p)}}\,dx <\infty.
\end{equation}
Note that this class depends on $p$ and $s$. The condition $s'=sp/2<1$ ensures that the class contains a large class of uniformly Lipschitz functions depending on our choice of $s$ and $p$. The value of $C(\var)$ only depends on the positivity, smoothness and  behaviour of $\var(x)$ as $|x|\to\infty$.

\medskip

\noindent{\bf Admissible decay rates.} It has been proved in \cite{BonfVaz2014MR3122168} that there are many smooth, bounded and positive functions $\varphi$ decaying at infinity like a power $\varphi\sim O(|x|^{-(N+\gamma)})$ with $\gamma>2s'$ such that ${\mathcal L}_{s'}\var$ decays like $O(|x|^{-(N+2s')})$ as $|x|\to \infty$, assuming that $0<s'<1$. We can check that  ${\mathcal M}_{s'}\var$ decays in the same way if $2s'<1$. We have
$$
\frac{|{\mathcal M}_{s'}\varphi(x)|^{1/(2-p)}}{\varphi(x)^{(p-1)/(2-p)}}   \sim |x|^{-\mu}
$$
with
$$
\mu=\frac{(N+2s')}{2-p}-\frac{(N+\gamma)(p-1)}{2-p}.
$$
The expression is integrable if $\mu>N$. Working out the details we find that $C(\varphi)$ is finite if
$\gamma<2s'/(p-1)$.  As convenient examples of admissible weights, we consider the family of weight functions
$$
\var(x)=(1+|x|^2)^{-(N+\lambda)/2}, \qquad \lambda>0.
$$
 Then we get a finite value for $C(\varphi)$ iff $\lambda<sp/(p-1)$. Actually, we can find  weights in the
 class $\mathcal C(s,p)$ with a decay at infinity that is closer to the limit power:
\begin{equation}\label{limit.growth.cond}
\var(x)\sim |x|^{-(N+(sp/(p-1))} (\log(1+|x|))^{\lambda}, \quad \mbox{with } \lambda>\frac{2-p}{p-1}.
\end{equation}
 We leave it to the reader to check this detail. The nonexistence result of Theorem \ref{thm.instbu} shows that the limit power cannot be reached.

\medskip

\noindent {\bf Numerical inequality.} We combine the simple inequalities valid for $a>b>0$ and $0<m<1$:
$$
a^m-b^m\le (a-b)^m, \quad a^m+b^m\le 2^{1-m}(a+b)^m,
$$
to conclude that for all possible values of real numbers $a$ and $b$ we have
\begin{equation}\label{est.numineq}
|a^m-b^m|\le 2^{1-m}|a-b|^m.
\end{equation}
Here  we  use the notation $a^m=|a|^{m-1}a$.

\subsection{Main estimate}
 We now state and prove our basic weighted estimate.

\begin{theorem}\label{them1} Let $0<s<1$, $1<p<2$, with $sp<1$. Let $u$ be a nonnegative semigroup solution of the FPLE \eqref{frplap.eq} in a strip $Q=\ren\times (0,T)$ with $T>0$. Then for all $\var\in \mathcal C(p,s)$  there is a finite constant $K>0$ depending only on $\var$ such that we have
\begin{equation}\label{est.wL1}
|X^{2-p}(t_1)-X^{2-p}(t_2)|\le K |t_1-t_2|.
\end{equation}
Actually, we may take $K(\varphi) =(2-p)C(\var)^{2-p}$ with $C(\var)$ given by \eqref{cond.phi}.
\end{theorem}

\noindent {\sl Proof.} (i) We  multiply an $L^2$ solution by a smooth and positive test function $\varphi(x)$, and integrate by parts in the equation as in \cite{Vazquez2020} to obtain for the evolution of the weighted mass \ :
$$
-\frac{d}{dt}\int u(x,t) \varphi(x) \,dx =\int {\mathcal L}_{s,p} u(x,t)  \,\varphi(x) \,dx=
$$
$$
\iint \frac{[u^{p-1}(x,t))-u^{p-1}(y,t)]\,(\varphi(x)-\varphi(y))}{|x-y|^{N+sp}}\,dxdy,
$$
where integrals extend to $\ren$. Therefore, by symmetry
$$
\big|\frac{dX}{dt}\big| \le 2\,\int u^{p-1}(x,t)\,dx\left(\int \frac{|\varphi(x)-\varphi(y)|}{|x-y|^{N+sp}}\,dy\right)=
2\,\int u^{p-1}(x,t)\, {\mathcal M}_{s'}\varphi(x)\,dx\,,
$$
where $s'=sp/2$. Since $m=p-1$ lies between 0 and 1, we get by H\"older's inequality
$$
\int u^{p-1}(x,t)\,{\mathcal M}_{sp/2}\varphi(x)\,dx\le \left(\int u(x,t)\varphi\right)^{p-1}dx\left(\int \frac{|{\mathcal M}_{sp/2}\varphi(x)\,dx|^{1/(2-p)}}{\varphi(x)^{(p-1)/(2-p)}}\,dx\right)^{2-p}\,.
$$
 In order to use this inequality we need to choose $\varphi$ so that this last expression, that we have called $C(\var)$, is finite. Taking $\var\in \mathcal C(p,s)$, we can write
$$
|\frac{dX}{dt}| \le X(t)^{p-1}\,C(\varphi)^{2-p}.
$$
Integration in time of this differential inequality gives the result   \eqref{est.wL1}. For semigroup solutions in other $L^q$ spaces use approximation.\qed

\medskip


\medskip
 The following result is a rescaled version of estimate \eqref{est.wL1}.

\begin{corollary}\label{cor.1} Under the same conditions, if $\var_R(x)=\var(x/R)$ we get
\begin{equation}\label{cond.phi.R}
\big|(\int u(t_2)\var_R\,dx)^{2-p}-(\int u(t_1)\var_R\,dx)^{2-p}\big|\le K(\varphi) R^{N(2-p)-sp}|t_1-t_2|.
\end{equation}
It follows that any integrable solution the finite mass is conserved in the range $p> p_c=2N/(N+s)$.
\end{corollary}

 \noindent {\sl Proof.} We only have to observe that for all $R>0$
$$
\mathcal L_{sp}(\varphi_R)(x)= R^{-sp}(\mathcal L_{sp} \varphi) (x/R),\quad
\mathcal M_{sp}(\varphi_R)(x)= R^{-sp}(\mathcal M_{sp} \varphi) (x/R),
$$
and use the formula for $K(\var)$. For the second part check that $p>p_c$ implies that $R^{N(2-p)-sp}\to 0$ as $R\to\infty$.\qed

\begin{corollary}\label{cor.2} Theorem \ref{them1} and Corollary \ref{cor.2} are true for differences of ordered solutions. If we consider two ordered solutions $0\le u_1\le u_2$, then estimate \eqref{cond.phi.R} holds for  $u=u_2-u_1$.
\end{corollary}

For the proof we only have to repeat the proof of the theorem after subtracting the two equations, and use the inequality $(a^m-b^m)\le (a-b)^m $ when $a>b>0$ and $0<m<1$. This means that $u_2^{p-1}(x,t)-u_1^{p-1}(x,y)\le u^{p-1}(x,y)$. The rest holds. \qed

\medskip

\noindent {\bf Remark on universality.}
It is very important that this expression can be applied to all nonnegative solutions (that can be suitably approximated, and  the constant $K$  does not depend at all on the solution. We say that it is a universal estimate. Note also that $t_1$ and $t_2$ are not required to be ordered.

\medskip

\noindent {\bf Comparative remarks.} (1)
 A similar but stronger  type of nonlinear estimate was found for the Fast Diffusion range of the Porous Medium Equation, cf. \cite{HerrPierre1984MR797051}. In that result $\var$ could be allowed to have compact support. This in turn allowed them to establish existence of solutions for the Cauchy problem in the whole space without any growth requirement in the initial data, it works for just locally integrable functions. Such generality is not true for our problem,  as we will show.

(2) We extended the local estimate to the fractional version of the PME in \cite{BonfVaz2014MR3122168}. There it lost its purely local form and took the weaker form on weighted integrals, that allowed for an existence theory with definite growth restrictions. Now weights with compact support are not allowed. In fact, unrestricted growth  is not expected in fractional diffusion, so the strong form cannot be true. The parameter restrictions were $0<s<1$ and $0<m<1$, no further restriction.

(3) On the other hand, in the case of the fractional linear heat equation the optimal growth rate for initial data in the Cauchy problem is just $\gamma<2s$, which will be in agreement with our results in the next section, derived from our basic estimate. The linear case was proved in \cite{BonSirVazMR3614666}. The proof relies on representation with a Green function that decays at infinity as expected, i.e., like $O(|x|^{-(N+2s)})$.

(4) No version of this inequality has been found  for the standard local model of the $p$-Laplacian evolution equation, $u_t=\mbox{div} (|\nabla u|^{p-2}\nabla u)$, and this fact complicated the study of the Cauchy Problem under optimal initial conditions done in \cite{DiBenHerrMR1066761}. It seems to us that the weighted $L^1$ estimate  fails.

\section{General existence theory}\label{sec.ex}

We keep the condition $sp<1$ in this section.
We are going to prove existence of solutions under quasi-optimal conditions on the initial data. The concept of weak solution is a function $u\in L^1(Q_T)$ where $Q_T=\ren\times (0,T)$ such that
\begin{equation}\label{def.ws}
\iint u\,\psi_t\,dxdt+\iiint \frac{(u(x,t)-u(y,t))^{p-1}(\psi(x,t)-\psi(y,t))}{|x-y|^{N+sp}}\,dxdydt=0
\end{equation}
for every smooth function $\psi\in C^\infty([0,T]\times\ren)$ with compact support in $Q_T$ so that in particular $\mathcal M_{sp/2}\psi$ is bounded. The integrals are extended to $\ren$ in $x$ and $y$ and to $(0,T)$ in time. The second term must be well defined and this will depend on the a priori estimates, as the proof below shows.

\begin{theorem}\label{thm.exist} {\rm [Existence]} Let $sp<1$ and let us consider the Cauchy Problem for the FPLE with locally integrable initial data $u_0\ge 0$. If
\begin{equation}\label{est.excond}
\int u_0(x)\,\var(x)\,dx<\infty
\end{equation}
for some admissible test function $\var\ge0$   in the class  $\mathcal C(s,p)$, then  there exists a weak solution of the problem which is defined in $Q_\infty$.  This solution is continuous in the weighted space, $u\in C([0,T]:L^1(\ren, \varphi\dx))$\,. The initial data are taken in the sense of strong convergence in $L^1_{loc}(\ren)$.
\end{theorem}

%

\noindent {\sl Proof.~} (i) We follow the outline of proof of Theorem 3.1 of \cite{BonfVaz2014MR3122168} for the fractional porous medium equation, but some important changes are needed. Let $\varphi \in \mathcal C$ and let $\varphi_R$ the scaling of $    \var$.
 Let $0\le u_{0,n}\in L^1(\ren)\cap L^\infty(\ren)$ be a non-decreasing sequence of initial data $u_{0,n-1}\le u_{0,n}$, converging monotonically to $u_0\in L^1(\ren, \varphi\dx)$. By the Monotone Convergence Theorem, it follows that $\int_{\ren}(u_0- u_{n,0})\varphi\dx \to 0$ as $n\to \infty$.

(ii) We prove existence of the monotone limit of the approximating solutions. Consider the unique strong solutions $u_n(t,x)$ of equation \eqref{frplap.eq} with initial data $u_{0,n}$, as constructed in \cite{VazFPL2-2020}. By the comparison results of that paper we know that the sequence of solutions is a monotone sequence. The weighted estimates of previous section imply that the sequence is bounded in $L^1(\ren, \varphi\dx)$ uniformly in $t\in[0,T]$\,.
 \begin{equation}\label{HP.unif}\begin{split}
&\left(\int_{\ren}u_n(t,x)\varphi(x)\dx\right)^{2-p}\le
 \left(\int_{\ren}u_n(0,x)\,\varphi(x)\dx\right)^{2-p} + K(\var)\,t\\
&\le \left(\int_{\ren}u_0(x)\varphi(x)\dx\right)^{2-p}
+ K(\var) \,t.
\end{split}
\end{equation}
  By the monotone convergence theorem in $L^1(\ren, \varphi\dx)$, we know that the solutions $u_n(t,x)$ converge monotonically as $n\to \infty$ to a function $u(t,x)\in L^\infty ((0,T): L^1(\ren, \varphi\dx))$. We also have
\begin{equation}\label{HP.s.3}
\left(\int_{\ren}u(t,x)\varphi(x)\dx\right)^{2-p}\le \left(\int_{\ren}u_0(x)\varphi(x)\dx\right)^{2-p}
+  K(\var)\,t
\end{equation}

(iii) We show next that the obtained limit function $u(t,x)$  is a weak solution to equation \eqref{frplap.eq} in $[0,T]\times \ren$ in the sense of Definition \ref{def.ws}.  We know that each $u_n$ is a bounded strong solution according to the theory of \cite{VazFPL2-2020}  since the initial data $u_0\in L^1(\ren)\cap L^\infty(\ren)$. Therefore, for all $\psi\in C_c^\infty([0,T]\times\ren)$ we have
\[\begin{split}
&\int_0^T\int_{\ren}u_n(t,x)\psi_t(t,x)\dx\dt=\\
&\int_0^T\iint \frac{(u_n(x,t))-u_n(y,t))^{p-1}(\psi(x,t)-\psi(y,t))}{|x-y|^{N+sp}}\,dxdydt.
\end{split}\]
Taking the limit $n\to\infty$ in the first line is easy:
\[
\lim_{n\to\infty}\int_0^T\int_{\ren}u_n(t,x)\psi_t(t,x)\dx=\int_0^T\int_{\ren}u(t,x)\psi_t(t,x)\dx\,,
\]
since $\psi$ is compactly supported and we already know that $u_n(t,x)\to u(t,x)$ in $L^1_{\rm loc}$.

(iii') On the other hand, the integral in the  second line, $I(u_n)$, is  well defined and can be estimated uniformly. We argue as before, by using the inequality
 $$
 |u_n(x,t)-u_n(y,t)|^{p-1}\le u_n(x,t)^{p-1}+u_n(y,t)^{p-1},
  $$
and  bounding the first of the two ensuing integrals by
$$
I_1(u_n)=\int dt \big(\int u_n(x,t)^{p-1}\,{M}(x,t)\,dx\big).
$$
where
$$
M(x,t)=\int \frac{|\psi(x,t)-\psi(y,t)|}{|x-y|^{N+sp}}\,dy.
$$
Due to the regularity of $\psi$, the last integral is bounded above by some $\overline M_\psi(x)$ that behaves like $C(1+|x|)^{N+sp}$ independently of $t$.
Hence,
$$
I_1(u_n)\le \int dt \left( \int u_n(x,t)\var\,dx\right)^{p-1} \, \left(\int \frac{ \overline M_\psi(x)^{1/(2-p)} }{\var^{(p-1)/(2-p)}}\, dx\right)^{2-p}
$$
so that finally
$$
I_1(u_n)\le C(\psi,\var)\int dt (\int  u_n(x,t)\var\,dx)^{p-1}  \le C\, T.
$$
%

The second  integral,  $I_2(u)$, is treated similarly  by exchanging $x$ and $y$.

(iii'') We  recall that  $u\ge u_n$, and we may apply the argument of the previous paragraph to $v_n=u-u_n$
as follows
$$
I(u)-I(u_n)=\iint \frac{(A(x,y,t)(\psi(x,t)-\psi(y,t))}{|x-y|^{N+sp}}\,dxdydt
$$
with
$$
A(x,y,t)=(u(x,t)-u(y,t))^{p-1}-(u_n(x,t)-u_n(y,t))^{p-1}.
$$

Using the numerical inequality \eqref{est.numineq} we conclude that for all possible values of $u(x,t))-u(y,t)$ and $u_n(x,t))-u_n(y,t)$ we have
$$
|A(x,y,t)|\le c(p)|u(x,t))-u(y,t)-(u_n(x,t))-u_n(y,t))|^{p-1}\le |v_n(x,t)-v_n(y,t)|^{p-1}
$$
and then we continue as before to prove that
$$
|I_1(u)-I_1(u_n)|\le C  \left(\iint  (u-u_n)\varphi\,dxdt\right)^{p-1}\left(
\iint \frac{|M |^{1/(2-p)}}{\var^{(p-1)/(2-p)}}\,dxdt\right)^{2-p}
$$
so that by virtue of the previous estimates on the weighted convergence of $u-u_N$ we get
$$
|I(u)-I(u_n)|\to 0 \quad \mbox{as \ } \quad n\to\infty.
$$

(iv) The solutions constructed above for $0\le u_0\in L^1(\ren, \varphi\dx)$ satisfy the weighted estimates \eqref{est.wL1} so that
\begin{equation}\label{HP.s.44}
\left|\int_{\ren}u(t,x)\varphi(x)\dx- \int_{\ren}u(\tau,x)\varphi(x)\dx\right|\le
2^{\frac{1}{1-m}}C_1 \,|t-\tau|^{\frac{1}{1-m}}
\end{equation}
which gives the continuity in $L^1(\ren, \varphi\dx)$\,.
Therefore, the initial trace of this solution is given by $u_0\in L^1(\ren, \varphi\dx)\,.$ \qed

\nc

\subsection{Uniqueness of limit solutions}
 We prove that the limit is independent of the approximating sequence.

\begin{theorem}\label{thm.uniqlarge} {\rm [Uniqueness]}
The solution constructed in Theorem \ref{thm.exist} by approximation from below is unique. We call it the minimal solution. In this class of solutions the standard comparison result holds, and also the weighted $L^1$ estimate of Theorem $\ref{them1}$\,. Extinction in finite time is not excluded.
\end{theorem}

\noindent {\sl Proof.~} We keep the notations of the proof of existence Theorem \ref{thm.exist}. We follow an argument that is well-known in the proof of uniqueness of minimal solutions of different problems, see \cite{GalVaz1997, GalVaz2002}, where they are also called proper solutions.

Assume that there exists another sequence $0\le v_{0,k}\in L^1(\ren)$ which is monotonically non-decreasing and converges monotonically to $u_0\in L^1(\ren, \varphi\dx)$\,. By the same considerations as in the proof of existence theorem, we can show that there exists a limit solution $v(t,x)\in C([0,T]:L^1(\ren, \varphi\dx))$. We want to show that $u=v$, where $u$ is the solution constructed in the same way from the sequence $u_{0,n}$. We will prove equality by proving that $v\le u$, and then the same argument will prove $u\le v$. To prove that $v\le u$ we use the estimates
\begin{equation}\label{final.001}
\int_{\ren}\big[v_k(t,x)-u_n(t,x)\big]_+\dx\le \int_{\ren}\big[v_k(0,x)-u_n(0,x)\big]_+\dx
\end{equation}
which hold for any $u_n(t,\cdot), v_k(t,\cdot)\in L^1(\ren)$, see \cite{VazFPL2-2020}. Letting $n\to \infty$ we get that
\[
\begin{split}
&\lim_{n\to\infty}\int_{\ren}\big[v_k(t,x)-u_n(t,x)\big]_+\dx
\le \lim_{n\to\infty}\int_{\ren}\big[v_k(0,x)-u_n(0,x)\big]_+\dx\\
&=\int_{\ren}\big[v_k(0,x)-u_0(x)\big]_+\dx=0,
\end{split}
\]
since $v_k(0,x)\le u_0$ by construction. Therefore also $v_k(t,x)\le u(t,x)$ for $t>0$, so that in the limit $k\to \infty$ we obtain $v(t,x)\le u(t,x)$\,. The inequality $u\le v$ can be obtained simply by switching the roles of $u_n$ and $v_k$\,.\qed

\begin{corollary}\label{cor.limitcond} There exists a unique minimal solution of the Cauchy problem for the FPLE for all locally integrable initial data $u_0\ge 0$ such that
\begin{equation}\label{est.ex.exam}
\int_{\ren}  u_0(x)\, \frac{\log(2+|x|)^{\lambda}}{(1+|x|)^{-(N+ \frac{sp}{p-1})}}\,dx<\infty
\end{equation}
for some $\lambda>\frac{2-p}{p-1}$. We can state a sufficient condition as an average growth rate: for all large $R>1$ we have
\begin{equation}\label{est.ex.exam2}
\oint_{B_R} u_0(x)dx:= R^{-N}\int_{B_R} u_0(x)dx\le C\,R^{sp/(p-1)}(\log R)^{-\mu},
\end{equation}
with $\mu>\lambda+1>1/(p-1)$. In particular, this holds if $u_0\in L^1_{loc}(\ren)$ and $u_0(x)\le C|x|^{\frac{sp}{p-1}}\log (1+|x|)^{-\mu}$ as $|x|\to\infty$.
\end{corollary}

\noindent {\bf Remark.} Positive  constants are examples of minimal weak solutions. This is easily deduced from the invariance of the set of approximating data under space translations.

\noindent {\bf Question.} We may wonder if it possible to find nonnegative solutions of the Cauchy Problem outside of this class of ``constructible minimal solutions''. The impression is that the answer is no, but such solutions exist in the theory of semilinear heat equations $u_t=\Delta u+ u^p$ for large values of $p>1$, as proved in \cite{GalVaz1997}, see also \cite{FilaMizomultiplecont2007}. They happen in blow-up situations (so-called incomplete blow-up).

\subsection{Singular minimal solutions}\label{sec-singms}

We will need to consider solutions in cases where the initial function is not locally integrable at some places with singular values. In the sequel we address the construction by the same method based on passing to the limit the approximation from below by  solutions with data in  the standard Lebesgue spaces. In that case we have near those singular places,
$$
\int_{B_R(x_0)} u_0(x)\,dx=\infty.
$$
Applying the weighted $L^1$ inequality with $\var=1$ in $B_R(0)$ to the approximations $u_n$, we have in the limit for every $t>0$
$$
\left(\int_{B_R(x_0)} \lim_n u_n(x,t)\var(x)dx\right)^{2-p}=\lim_n C_n -K(\var)t^{}=\infty.
$$
hence the limit solution is not locally integrable at those place. We will see  examples of such singular minimal solutions where we prove that the limit is finite away from the space singularity, and examples where the limit is infinite everywhere for $t>0$ (what we will call instantaneous blow-up). See more in Subsection  \ref{sec.sing.sss}.

\subsection{Solutions with changing sign}\label{sec.changsing}

Though our main interest lies with nonnegative solutions, we add the information on existence of weak solutions
for data and solutions with changing sign. We arrive at the idea of min-max and max-min solutions. The main estimate still holds for ordered differences.

\begin{theorem}\label{thm.exist.signed} Let us consider the Cauchy Problem for the FPLE with locally integrable signed initial data $u_0$. If
\begin{equation}\label{est.excond.sgn}
\int |u_0(x)|\,\var(x)\,dx<\infty
\end{equation}
for some admissible test function $\var\ge0$   in the class  $\mathcal C(s,p)$, then  there exists a weak solution of the problem which us defined in $Q_\infty$.  This solution is continuous in the weighted space, $u\in C([0,T]:L^1(\ren, \varphi\dx))$\,. The initial data are taken in the sense of strong convergence in $L^1_{loc}(\ren)$.
\end{theorem}

\noindent {\sl Proof.~} (i) We follow the method of construction by approximation. First, we observe that both
$f=\max\{u_0\}$, and $g=\max\{-u_0,0\}$ satisfy the hypothesis of Theorem \ref{thm.exist}, so there are monotone sequences of good data (good in the sense of the proof of the theorem) such that $0\le f_n(x)$ converges to $f$, $0\le g_m$ converges to $g$, and the corresponding solutions $U_n(x,t)$ and $V_m(x,t)$ converge to the respective minimal weak solutions $W(x,t)$ and $V(x,t)$. Note that $u_0=f-g$.

(ii) Now we consider the simpler case where $g$ is an $ L^1$ function. In that case we do not need approximation from below, just put $g_m=g$ and $u_{m,n} = -g+ f_n(x),$ where $f_n$ is a monotone approximation to $u_0+g=f$.
then we may repeat the proof of Theorem \ref{thm.exist} with minor changes and obtain a minimal weak solution, that will be unique in the sense of Theorem  \ref{thm.uniqlarge}.
A similar argument solves the case where $f$ is an $ L^1$ function, and then we obtain a unique maximal weak solution by downwards approximation. These are just extensions of Theorem \ref{thm.exist}.

(iii) Let us tackle the general case. We consider a double sequence of approximations $u_{nm}(x,0)$ such that $-g_m(x)\le u_{n,m}(x,0)\le f_n(x)$, the sequence is increasing in $n$ and decreasing in $m$, and
$$
\lim_{n\to\infty} u_{m,n}(x,0)=f(x), \quad \lim_{n\to\infty} u_{m,n}(x,0)=-g(x).
$$
We examine the corresponding solutions $u_{m,n}(x,t)$ and their limits. We easily obtain  $-V_m(x,t)\le u_{m,n}(x,t)\le W_n(x,t)$,
which is a comparison between good and minimal solutions. Therefore, in the limit $n\to\infty$ we obtain a function $U_1(x,t)$ such that
$$
-V_m(x,t)\le u_{m,n_1}(x,t)\le \lim_{n\to\infty} u_{m,n}(x,t)=U_{1,m}(x,t)\le W(x,t),
$$
This $U_1$ is the large solution in the sense of point (ii) to initial data $U_{1,m}(0)$, we can check that $U_{1,m}(0)\le f-g$. Since the sequence is monotone downwards in $m$ we get a further limit
$$
U_1(x,t)=\lim_{m\to\infty}\lim_{n\to\infty}u_{m,n}(x,t)=U_1(x,t),
$$
and we easily check that $-V(x,t)\le U_1(x,t)\le W(x,t)$. This is the candidate to max-min-solution.
It is a weak solution and takes on the initial data as in Theorem \ref{thm.exist}.\nc

(iv) Reversing the order of the limits we obtain
$$
-V(x,t)\le  \lim_{m\to\infty} u_{m,n}(x,t)=U_{2,m}(x,t)\le u_{m_1,n}(x,t) \le W(x,t),
$$
and
$$
U_2(x,t)=\lim_{n\to\infty}\lim_{m\to\infty}u_{m,n}(x,t)=U_2(x,t)
$$
and we easily check that
$$
-V(x,t)\le U_2(x,t)\le U_1(x,t)\le W(x,t).
$$
This $U_2$ is the candidate to min-max-solution.\qed

\medskip

\noindent {\bf Question.} We do not have a proof of uniqueness, $U_1=U_2$, unless there is a good bound on the initial data from below or above. Can we have non-unique solutions?

\noindent {\bf Example.} There is a typical example of signed solution, i.e., the linear function, $u_c(x,t)=c_1x_1$ (and its translates and rotations). Note  that it is constant in time. In case $sp/(p-1)>1$ this solution comes under the scope of Theorem \ref{thm.exist.signed}, but for $sp\le p-1$ it does not. For $s<1/2$ both possibilities arise depending on $p$. However, for $s\ge 1/2$ we always have the compatibility $1< sp/(p-1)$.

\noindent  {\bf Question.} Is the linear solution, constant in time but not in space, reachable by the above approximation process? We know that for $sp\le p-1$ it is not. \nc

\section{Positivity versus extinction}\label{sec.masscon}

We will produce in the sequel a number of solutions, even semigroup solutions, that vanish in finite time, like the example \eqref{decay.vss.u2}, and we will devote a section to discuss the issue. On the other hand, (strict) positivity is usually a given property of nonnegative and nontrivial solutions of heat equations and fast diffusion equations, while the property is lost for slow diffusion equations because of the property of finite propagation and the existence of free boundaries, \cite{FriedBkMR679313}.

Here, we want to settle the dichotomy between positivity and extinction for  nonnegative solutions of the FPLE. For any solution $u(x,t)\ge 0$ and  at any given time $t>0$, the alternative that we describe in Theorem \ref{thm.mc} holds for $u(\cdot, t)$ as a function of $x$. The theorem  allows us to define the extinction time $T(u_0)$  as the first time where $u(x,t_1)$ is the trivial function, and equivalently, as the first time where the nonnegative solution is no more strictly positive and touches zero.

Here, $sp<1$ is not required. Fistly, we need a technical lemma on mass control that has other uses.

\begin{lemma}\label{lemma.mc} \rm [Conditional Mass conservation] \sl  Let $u(x,t)\ge 0$ be the semigroup solution of the Cauchy Problem with initial data $u_0\in L^1(\ren)$, $u_0\ge 0$, and assume that $u(x,t_1)$ is compactly supported for some $t_1>0$. Then the mass is conserved for all $t\ge t_1$.
\end{lemma}

\noindent {\bf Remark.} Conservation of mass was proved in \cite{VazFPL2-2020} without the compact support assumption in the range $p_c<p<1$. In view of the examples of finite-time extinction, \bc the result does not hold when $1<p<p_c$. \nc

\medskip

\noindent {\sl Proof of the lemma.} (i) {\sc Reduction step.} 
We may always assume that $u_0\in L^1(\ren)\cap L^\infty(\ren)$ and $u_0$ is compactly supported. If our form of mass conservation is proved under these assumptions, then it follows for all data $u_0\in L^1(\ren)$ by the semigroup contraction property.

Let $B_R(0)$ contain the support of $u(\cdot,t_1)$. Due to the pointwise inequality $u_t\le u/(p-2)t$, we get the conclusion that whenever $u(x_1,t_1)=0$, then we have $u(x_1,t)= 0$ for all $t>t_1$ so that the vanishing set is preserved in forward time.  Therefore,  the support of $u(\cdot, t)$ will be contained in the same ball $B_R(0)$ for all $t\ge t_1$.

(ii) Let us assume that $sp<N$, something that will always happen for $N\ge 2$.
We do a calculation for  the tested mass. Taking a smooth and compactly supported test function $\varphi(x)\ge0 $, we have for $t_2>t_1>0$:
\begin{equation}\label{mass.calc}
\left\{\begin{array}{l}
\displaystyle \left|\int u(t_1)\varphi\,dx-u(t_2)\varphi\,dx\right|\le \iiint
\left|\frac{\Phi(u(y,t)-u(x,t))(\varphi(y)-\varphi(x)}{|x-y|^{N+sp}}\right|\,dydxdt\\[10pt]
\le \displaystyle \left(\iiint  |u(y,t)-u(x,t)|^p\,d\mu(x,y)dt\right)^{\frac{p-1}{p}}
\left(\iiint |\varphi(y)-\varphi(x)|^p\,d\mu(x,y)dt\right)^{\frac{1}{p}}\,,
\end{array}\right.
\end{equation}
with space integrals over $\ren$  and time integrals over $[t_1, t_2]$. Use now the sequence of test functions $\varphi_n(x)=\varphi(x/n)$ where $\varphi(x)$ is a cutoff function which equals 1 for $|x|\le 2$ and zero for $|x|\ge 3$. We take $n\ge R$. Then, have to consider different regions for the calculation with the multiple integrals. Note that in \eqref{mass.calc} we estimate integrals in absolute value (by taking absolute value of the integrand)\nc.

 We first deal with exterior region $A_R=\{(x,y): |x|,|y|\ge R\}$. Recalling \eqref{mass.calc} and the assumption on the compact support made on $u$ we have
\begin{equation*}\label{mass.calc2}
 \displaystyle I(A_n):= \int_{t_1}^{t_2}\iint_{A_n} \frac{|\Phi(u(y,t)-u(x,t))|\,|\varphi_n(y)-\varphi_n(x)|}{|x-y|^{N+sp}}\,dydx\,dt=0.
 \end{equation*}
Arguing in a similar fashion in the inner region $B_n=\{(x,y): |x|,|y|\le 2n\}$ where
$\varphi_n(x)-\varphi_n(y)=0$, we see that its contribution to the integral \eqref{mass.calc} is also zero.

 We still have to make the analysis in  other regions so that we cover the whole domain $x,y\in \ren$. An option of to consider the cross regions
$C_n=\{(x,y): |x|\ge 2n ,|y|\le R\}$ and $D_n=\{(x,y): |x|\le R ,|y|\ge 2n\}$. Both are similar so we will look only with $D_n$. The idea is that we have an extra estimate: $|x-y|>n$ that avoids the singularity in the weight of the integrand. We have
\begin{equation*}
 \begin{array}{c}
  \displaystyle I(D_n)\le \int_{t_1}^{t_2}\iint_{D_n} |u(x,t)|^{p-1}(1-\var_n(y))\,d\mu(x,t)dt\le \\
\displaystyle \int_{t_1}^{t_2}dt \big(\int_{B_R} dx \,|u(x,t)|^{p-1}
  \big(\int_{|x-y|>n} |x-y|^{-N-sp} dy\big)   \big)\\
 \displaystyle  \le Cn^{-sp}\int dt \int_{B_R} |u(x,t)|^{p-1}\,dx)\,.
 \end{array}
\end{equation*}
Since $0<p-1<1$, we have
$$
\int_{B_R} |u(x,t)|^{p-1}dx\le R^{N(2-p)}\left(\int_{B_R} |u(x,t)|\,dx\right)^{p-1}
\le R^{N(2-p)}\|u(x,t_1)\|^{p-1}.
$$
Therefore, $ I(D_n)\le K n^{-sp} $ which tends to zero as $n\to\infty$ with a power rate. Same  for $I(C_n)$. This concludes the proof. Note that these regions overlap but that is no problem.

(iii) The case $N=1$ and $sp\ge 1$. Conservation of mass in then quite easy, see a proof in \cite{VazFPL2-2020}. \qed

\begin{theorem}\label{thm.mc} Let $1<p<2$.  Let $u(x,t)\ge 0$ be the semigroup solution of Problem with initial data $u_0\in L^1(\ren)$, $u_0\ge 0$. Then for every $t_1>0$  the space function $u(\cdot, t_1) $ is either strictly positive everywhere or identically zero.
\end{theorem}

\medskip

\noindent {\sl Proof.} We will use the continuity of the semigroup map in $L^1(\ren)$ and its monotonicity to propose a first extinction time definition as
$$
T_1(u_0)=\sup\{t>0: \|u(t)\|_1>0\}.
$$
Of course, it may happen that $T_1(u_0)=\infty$. This is the case for all $u_0\in L^1(\ren)$ if $p>p_c$ since we already know that there is mass conservation. But there are many examples of initial data with  finite extinction time, $T_1(u_0)=\infty$, if $p<p_c$.
$T_1$ cannot be zero by the continuity in time of the norm $\|u(t)\|_1$,  proved in the construction of the semigroup.

(i) Let us assume moreover that $u_0$ is bounded and compactly supported, say in the ball of radius $R>0$. Then, by using Aleksandrov's Maximum Principle, we know that for every $t>0$ we have monotonicity in radial outward cones of directions of space  with vertex at any point \ $|x_1|\ge 2R$,  and also \ $u(x_2,t)\le u(x_1,t)$ if $|x_2|\ge |x_1|+2R$. Therefore, if at some point $x_1$ with $|x_1|\ge 2R$  we have $u(x_1,t_1)=0$, then we also have $u(x,t)=0$ for all $ t\ge t_1$ and $|x|\ge 4R$. In this situation, the previous lemma guarantees conservation of mass for times $t\ge t_1$.

(ii) Assume now that  $T_1(u_0)$ is finite. By the continuity of the mass and its constant value for $t\ge t_1$, it follows that there can be no extinction at a later time unless $u(x,t_1)$ is already the trivial solution, hence $T_1=t_1$.  In that case the argument shows that uniform positivity is ensured at $t=t_1$ for all $C\ge |x|\ge 2R$ for any large $C$. Using the partial monotonicity of the solution in time $(2-p)tu_t\le u$, and in space (along outward cones, as explained above) we conclude that the solution is positive, actually uniformly positive locally in compact subsets  of the region
$$
\{(x,t): |x|\ge 2R, \ 0<t<T_1\}.
$$

(iii) We also have to consider the possibility $T_1=\infty$. If that would be the case, the assumption $u(x_1,t_1)=0$ leads to a solution supported in the ball of radius $4R$ in the interval $t_1\le t<\infty$. We have just proved that in that case the mass would be conserved. But general theory proves that $\|u(t)\|_2$ must go to zero for such a solution, \nc and this means (by compact support) that also $\|u(t)\|_1\to 0$.

Moreover, this contradiction leads to the positivity in the same outer sets as before, where now the time interval is $0<t<\infty$.

(iv) We now address the positivity in the initial core $B_R$. Let now $0< t_1<T_1(u_0)$ (finite or infinite) and consider the solution at $t_1$ which is positive around a point $x_1$ such that $|x|_1>2R$. There is a small constant and a small radius such that $u(x,t_1)\ge c$ in $B_r(x_1)$.

We now consider the semigroup solution $u_1$ with such lower data $u_1(x,t_1)=c\chi_{B_r(x_1)}(x)$. It  will have the form $u_1=F(|x-x_1|,t-t_1)$  with $F$ radially decreasing in space. Applying the argument of the previous step, there will a small time $\ve$ and a constant $c_1$ such around any other point $x_2\ne x_1$ we have
$$
u(x,t+\ve)\ge u_1(x,t+\ve)= F(x-x_1,\ve)\ge c_1
$$
We have concluded that $u(x,t)$ is strictly positive locally in all of $\ren$ for  $t_1+\ve>0$. This happens also  for smaller times (away from  0) by partial monotonicity. But $t_1+\ve$ is as close to $T_1$ as we want.

(v) We eliminate the assumption of boundedness and compact support on $u_0$ by monotone approximation from below.

\medskip

\noindent {\bf Extinction definition.} In view of these results we can propose equivalent definitions of Extinction Time for nonnegative solutions of the FPME in the whole space:
$$
T(u_0)=\max\{ t>0: \ u(x,t) \ \mbox{is positive for all} \ x\in\ren\}.
$$
or
$$
T(u_0)=\min\{ t>0: \ u(x,t)= 0 \ \mbox{for all} \ x\in\ren\}.
$$
Both definitions are equivalent to $T_1(u_0)$ defined above.
\medskip

\noindent {\bf Remark on positivity in other scenarios.} Positivity of nonnegative, nontrivial solutions is well-known in the linear fractional case $p=2$ thanks to the representation theorem, cf. \cite{BonSirVazMR3614666}. It also true for the superlinear case $p>2$, though the proof is less immediate, cf. \cite{Vazquez2020}. In both cases mass is conserved and no extinction in finite time may occur.

Positivity versus extinction occurs also for the Fractional Fast Porous Medium Equation and our present proof applies to that case. It also applies  to the standard non-fractional $p$-Laplacian equations for $p< 2$. The result is false for $u_t=\mbox{\rm div} (|\nabla  u|^{p-2}\nabla u)$ with $p>2$ due to the phenomenon of finite propagation and the presence of free boundaries, see for instance \cite{KaminVaz1988}.

\section{Scaling, power-like functions and self-similarity}\label{sec.self}

The arguments of this section are valid for $0<s<1$ and all $p>1$ unless a more restrictive range is indicated.
We begin here the study of self-similarity that will lead to the construction of a large set of new examples of special solutions of the FPLE equation \eqref{frplap.eq}. Already known examples are the constant functions, the fundamental solutions
$$
U(x,t)=t^{-N\beta}F(xt^{-\beta}), \quad \beta=1/(sp+N(2-p),
$$
constructed in \cite{VazFPL2-2020} in the optimal range $p_c<p<\infty$, as well as the very singular solutions (VSS) that we recall at the end of the section.

The ideas of self-similarity for nonlinear equations, like the porous medium equation, the $p$-Laplacian equation, and many other equations, as well the importance of the self-similar solutions to describe the asymptotic behaviour of more general solutions, are now very popular tools. Both topics were disseminated in Barenblatt's books like \cite{BarentSSSIA, BarenScalingBk}. The presentation below follows the ideas of the book \cite{Vazquez2006}.

\subsection{Scaling}\label{ssec.scaling}

In our study we will use the fact that the equation admits a scaling group that conserves the set of solutions, as mentioned in papers \cite{Vazquez2020, VazFPL2-2020}. Thus, if $u$ is a  strong solution of the equation, then we obtain a two-parameter family of solutions of the same type,
$$
\widehat u(x,t)= A u(Bx,Ct), \quad A, B, C>0,
$$
on the condition that $A^{2-p}C=B^{sp}$. This applies both to the semigroup solutions with data in the Lebesgue spaces constructed in \cite{Vazquez2020, VazFPL2-2020}, and to the minimal solutions in the large weighted spaces constructed above. See special choices of the scaling parameters in \cite{Vazquez2020, VazFPL2-2020}.
%
%
Let us also remark that the set of solutions of the equation is invariant under a number of isometric transformations, like: change of sign: $u(x,t)$ into $-u(x,t)$, rotations and translations in the space variable, translations in time, as well as vertical translations. They will also be used in the sequel.

\subsection{Scaling for power data}\label{ssec.scaling.prd}

 We want to investigate the existence of solutions of the FPLE \eqref{frplap.eq} with power-like data
\begin{equation}\label{ssdata}
u_0(x)=f(\theta)|x|^{\gamma}, \quad \theta=\frac{x}{|x|}\in \mathbb{S}^{N-1},
\end{equation}
where $\gamma$ can have both signs. The answer will depend of course on $\gamma$. We want to keep the invariance condition, $A^{2-p}C=B^{sp}$. If we also want to preserve the scaling of the power of such data we have $\widehat u(x,0)=u_0(x)$, i.e., $AB^{\gamma}=1$. This leads to the choices
$$
A= B^{-\gamma}, \quad C=B^{sp-\gamma (p-2)},
$$
with $B$ a free parameter. If we have a uniqueness theorem, as the ones mentioned above, we will have
$u(x,t)=\widehat u(x,t)$, and thus we obtain the  formula
$$
u(x,t)=B^{-\gamma} u(Bx, B^{sp-\gamma (p-2)}t) \qquad \mbox{for all } \ B>0.
$$

Suppose now that $sp-\gamma (p-2)\ne 0$.  There are two options.

\medskip

\noindent {\bf Self-similarity of the first kind.} If $sp+\gamma (2-p)>0$, which is  true when $1<p<2$ for
$$
\gamma> -\frac{sp}{2-p},
$$
i.e., either $\gamma\ge 0$, or  if the power is negative, it must be bounded below by an estimate that is at least $-s$. We then put
$$
sp+\gamma (2-p) =1/\beta.
 $$
 Then we obtain self-similarity with the usual trick: fixing $t$ and choosing $B$ so that  $t=B^{-1/\beta}>0$ we get
$$
u(x, B^{-1/\beta})= B^{-\gamma} u(Bx, 1).
$$
Putting now $B^{-1/\beta}=t$ again and $F(y)=u(y,1)$ we get
\begin{equation}\label{fss1k}
u(x, t)= t^{\gamma\beta} F(xt^{-\beta}).
\end{equation}
This is a form called self-similar solution of the first kind. The profile $F$ must satisfy
\begin{equation}\label{stateq.type1}
{\mathcal L}_{s,p}F(y)=-\beta\gamma\, F(y)+\beta y\cdot \nabla F(y).
\end{equation}
The problem is then reduced to find a suitable profile $F$. Once $F$ is found, the solution  $u$ will exist for all positive times. We may also insert an innocent time displacement
and write
$$
u(x, t; T)= (t-T)^{\gamma\beta} F(x\,(T-t)^{-\beta}).
$$
This solution will exist forward in time, starting at $t_0=T$. This similarity is also called forward self-similarity or standard self-similarity.

We can also have self-similarity of the first kind when $sp+\gamma (2-p)<0$, but then $\beta<0$ and the profiles tend to shrink with time according to formula \eqref{stateq.type1}. Such a phenomenon has been study if fast diffusion for the PME in \cite{Vazquez2006}.
We will concentrate here on the cases $\beta>0$ that we call {\sl expanding self-similar solutions of the first kind}. See more in Sections \ref{sec.exss} and \ref{sec.exss2}.

\medskip

\noindent {\bf Self-similarity of the second kind.}  If on the contrary $sp+\gamma (2-p)<0$ (which means $\gamma$ very negative), we may put
$$
\gamma (p-2)-sp=1/\beta'.
$$
Then, arguing as before and using inverse time, we get
\begin{equation}
u(x, t)= (T-t)^{-\gamma\beta'} F(x\,(T-t)^{\beta'}).
\end{equation}
This is a self-similar solution of the second kind, and applies in backward time $-\infty<t<T$, the constant $T$ being arbitrarily chosen. The profile $F$ must satisfy
\begin{equation}
{\mathcal L}_{s,p}F(y)=-\beta'\gamma F(y)+\beta' y\cdot \nabla F(y).
\end{equation}

\noindent {\bf Self-similarity of the third kind.} In the special case $\gamma_1= -sp/(2-p)$ (and $p\ne 2$), the scaling formula simplifies to
$$
u(x,t)=B^{-\gamma_1} u(Bx, t),
$$
which is just a scaling that preserves time.  So this way to self-similarity is barred. Note that for $1<p<2$ we have $\gamma<0$.

In this case there is another tool available, i.e., the self-similarity with exponential time factors.  We use the exponential Ansatz
\begin{equation}
u(x, t)= e^{- \alpha \,t} F(x\,e^{-\beta t}),
\end{equation}
Arguing as before, we eliminate time if the parameter $\beta$ is arbitrary and $\alpha= |\gamma_1|\beta$. This is a self-similar solution of the third kind, and applies for all times $-\infty<t<\infty$. The profile $F$ must satisfy
\begin{equation}
{\mathcal L}_{s,p}F(y)= \beta|\gamma_1| F(y)+\beta y\cdot \nabla F(y).
\end{equation}
We see that to find $F$ for $\beta=1$ is enough (by scaling, as pointed above). If we find $F$ we will obtain an eternal solution of the FPLE \eqref{frplap.eq}. Eternal solutions have been obtained for  local equations, $s=1$.

\medskip

\noindent {\bf Self-similarity in separated variables.}  There is another popular form of special solutions, called separate-variable solutions, of the form
$$
u(x,t)=A(t)F(x).
$$
We have constructed in \cite{VazFPL2-2020} solutions of this form that will be important for our studies below.
Note that the separate variable form is a solution of our equation is two separate equations hold: $A'(t)=-cA(t)^{p-1}$ and $\mathcal L_{s,p}F(x)=cF$.
The first one gives
$$
A(t)= at^{-1/(2-p)}
$$
while the second is a nonlinear eigenvalue problem that has been solved in \cite{VazFPL2-2020} in the form
given in its  Theorems 10.1 and 16.2 that we quote together here.

\begin{theorem}\label{thm.vss} (1) Let $p_c<p<p_1$. There exists a constant $C_\infty(s,p,N) >0$ such that
\begin{equation}\label{decay.vss.u1}
U_\infty(x,t)=\,C_\infty  \,(t+T)^{1/(2-p)} { |x|^{-sp/(2-p)} }
\end{equation}
is  a classical solution of equation \eqref{frplap.eq} for $x\ne 0$ which has a non-integrable singularity at $r=0$.  It is called the Very Singular Solution, VSS.

(2) Let $1<p<p_c$. There exists a constant $C_\infty(s,p,N) >0$ such that for any $T>0$ the function
\begin{equation}\label{decay.vss.u2}
U(x,t)=\,C_\infty  \,(T-t)^{1/(2-p)} { |x|^{-sp/(2-p)} }
\end{equation}
is a weak solution of the FPLE at all times $0<t<T$ and points $x\ne 0$ with an integrable singularity at $r=0$.
$U(x,0)\le L^q_{loc}(\ren)$ for all $1\le q< q_*=N(2-p)/sp$.
\end{theorem}

The solution in the first option lives for infinite time (it has forward similarity), the last one is an example of backwards similarity that vanishes in finite time.

These  will play a role in the sequel. In the case (1) they are called VSS and they are self-similar functions whose data are limit cases of integrable self-similar solutions.
\nc

\begin{figure}[t!]
 \centering
 \vspace{-2cm}
 \includegraphics[scale=0.45]{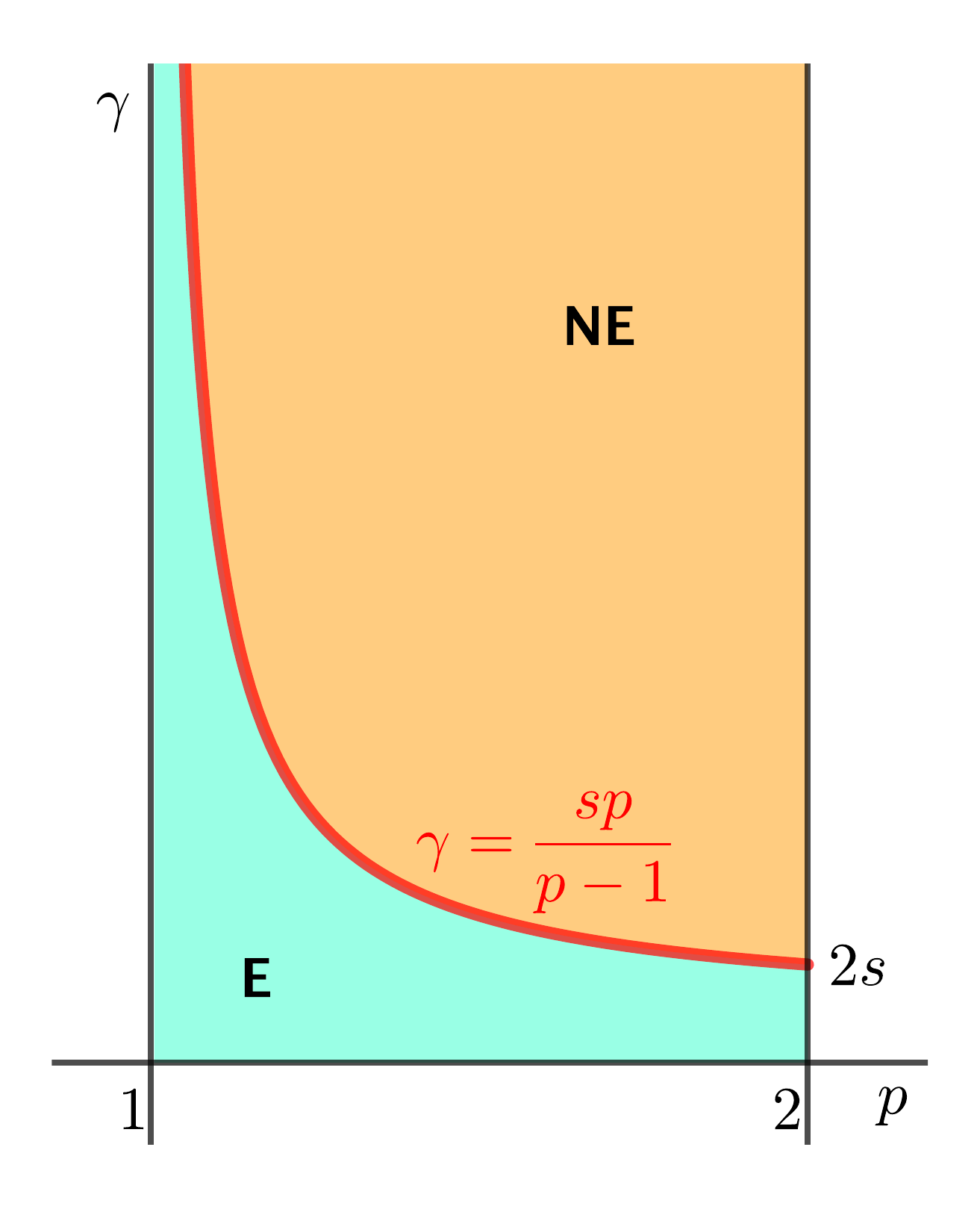}
 \caption{{Exponents of growing self-similar solutions. Region E represents existence, NE  denotes nonexistence, the separating line is also NE.}}
  \label{fig:2}
\end{figure}

\section{Self-similar solutions with growing data}\label{sec.exss}

Here we study the existence of self-similar solutions for growing data. We assume  for the moment that $sp<1$. We call $\gamma_1 =sp/(p-1)>0$, the growth exponent mentioned in Section \ref{sec.ex}. Note that $sp+\gamma(2-p)>0$ always since $\gamma\ge 0$.

\begin{theorem}\label{thm.ss.growing}  We consider the Cauchy problem for the FPLE \eqref{frplap.eq} with initial data
\begin{equation}
u_0(x)=A\,|x|^{\gamma},  \quad A,\gamma>0.
\end{equation}

\noindent (i) For every $\gamma\in (0,\gamma_1)$ with $\gamma_1=sp/(p-1)$ and every $A>0$, there exists a unique  minimal solution of this problem. It has the  self-similar form
\begin{equation}\label{ext.sst1}
u(x, t)= t^{\gamma\beta} F(xt^{-\beta})
\end{equation}
with $\beta=1/(sp+\gamma(2-p))>0$. The range of exponents $\gamma>0$ is optimal.

\noindent (ii) The case $A>0$ can be reduced to the case $A=1$ by the scaling formulas. If we denote by $u_A(x,t)$ the solution with initial data $A\,|x|^{\gamma}$, then
\begin{equation}
u_A(x,t)= A\,u_1(x,A^{p-2}t),    \qquad F_A(y)=A^{sp\beta}F_1(A^{(2-p)\beta}y).
\end{equation}
\noindent (iii) The profile $F$ is radially symmetric, continuous and positive with $F(0)>0$.
$F(x)$ is monotone increasing in $r=|x|$. As $|x|\to \infty$ the profile satisfies $F(y)\sim A|x|^{\gamma}$.
$F$ is a weak solution of the stationary equation \eqref{stateq.type1}.

\noindent (iv) We get \ $0\le F'(r)\le \gamma F(r)/r$, so that $F'=O(r^{\gamma-1})$ for large $r=|x|$.

\noindent (v) The solution $u$ is increasing in time, $u_t\ge0$, and $\mathcal L_{sp}(F)\le 0$.

\noindent (vi) The solution $u$ blows up in infinite time with rate
$$
u(x,t)\ge u(0,t)=F(0)\,t^{\alpha}, \qquad \alpha=\gamma\beta=\frac{\gamma}{sp+\gamma(2-p)}\,,
$$
for all $x\in \ren$ and $t>0$. Note that $\alpha$ increases with $\gamma$ from $\alpha(0)=0$ to
$\alpha(\gamma_1)=1$.
\end{theorem}

\noindent {\bf Remarks.} 1) The results are extended to the case $\gamma=0$ where we have the constant solution
$u(x,t)=A$ and the previous algebraic calculations agree.

2) Note that for the linear case $p=2$ we get  $\gamma_1=2s$, $\beta= 1/(2s)$ and $\alpha= \gamma/(2s)$. All of them are optimal values as proved in \cite{BonSirVazMR3614666}.

3) When $\gamma\to \gamma_1 $ we get $\beta(\gamma_1)=1/\gamma_1$, hence $\alpha(\gamma_1)=1$.

4) Moreover, as $p\to 1$  the $\gamma $-interval stretches to $[0,\infty)$, and we get in the limit $p\to1$ the exponents $\beta\to 1/(s+\gamma)$ and  $\alpha=\gamma/(s+\gamma).$

\medskip

\noindent {\sl Proof of the Theorem.~} We divide the proof in several steps, but we do not follow the same order as in the statement of the Theorem.

(i) If $sp<1$ the existence and uniqueness of a minimal weak solution is a consequence of Section \ref{sec.ex} and the self-similarity of the form \eqref{ext.sst1}
is a consequence of the analysis of Section \ref{sec.self}. Radial symmetry immediately follows. Monotonicity in $r=|x|$ comes from Aleksandrov symmetry principle. The weak formulation for $F$ also follows.

(ii) Because of Lemma  \ref{lemm.plcalc} we know that $u_t(0)\ge 0$ and then we will have $u(x,t)\ge u(x,0)$ for small $t$ and maximum principle will imply $u_t\ge 0$ for all $t>0$.
This is rigorously proved as a consequence of Lemma \ref{lemma.tech}. \nc This also implies that $u(x,t)\ge u_0(x)$ for all $x$, hence $F(y)\ge A\,|y|^{\gamma}$. It follows that $F$ is positive unless maybe at the origin.

(iii) We deduce  from $u_t\ge 0$ that $\mathcal L_{s,p}F(y)\le 0$. Using equation \eqref{stateq.type1}, we have $\gamma F(y)\ge yF'(y)$ so that (writing $r=|x|$)
the function $G(r)=F(r)r^{-\gamma}$ is nonincreasing, hence it has a limit at infinity that must be
$$
\lim_{r\to\infty}F(r)r^{\gamma}=A_1\ge A.
$$
 The fact that $A_1=A$ happens because it $u(x,t)$ takes the initial data according to the conclusions of Theorem \eqref{thm.exist}.

(iv) We need to check that $F(0)>0$ and that $F(r)$ is continuous at $r=0$.
For small $0<\gamma\le 1$ we argue as follows: the data are $\gamma$-H\"older with constant $A$. By the property of $L^\infty$ contraction (see \cite{VazFPL2-2020}), we deduce that so happens to $u(t)$, in particular to $F$ (take $t=1$), therefore we get
$$
F(r)-F(0)\le A\,r^\gamma.
$$
If now $F(0)=0$, the function $F$ would equal $u_0$,  and $u(x,t)$ be a stationary solution,so that $\mathcal L_{s,p}F=0$, which is not true. This means that $F(0)> 0$. It also follows from the formula that $F$ is continuous even at $r=0$.

To prove that $F_{\gamma'}(0)>0$ for  ${\gamma'}$ if $F_\gamma(0)>0$ for a certain $0<\gamma<{\gamma'}$, we just observe that at the initial time
$$
|x|^{\gamma'}\ge |x|^{\gamma}-C.
$$
for some $C$, and then we can compare  the solutions and conclude that in particular for $x=0$ we have
$$
u_{\gamma'}(0,t)=F_{\gamma'}(0)t^{\alpha}\ge  u_{\gamma}(0,t)-C=F_{\gamma}(0)t^{\alpha}-C.
$$
It follows that $F_{\gamma'}(0)\ge F_{\gamma}(0)>0$. Actually, we have proved that $F_\gamma(0)$ is monotone nondecreasing in $\gamma$.

(vi) A different argument to prove that $F$ is continuous at $r=0+$  uses a comparison argument based on displacement and is good for all $\gamma\in (0,\gamma_1)$. Let $A=1$,  and let $F(0)=a>0$ and $F(0+)=b>a$. Let us compare $u_1(x,t) $ with initial data $|x|^{\gamma}$ with $\overline u(x,t)$ initial data $\overline u(x,0)=(1-\ve)|x-x_0|^{\gamma}-C$, where $|x_0|=h$ is small and such that
$$
|x|^{\gamma}\ge (1-\ve)|x-x_0|^{\gamma}-C.
$$
Then $C$ is small if $h$ is small. We then have
$$
u_1(0,t)\ge \overline u(0,t)=u_{1-\ve}(h,t)-C= (1-\ve)u_1(h, (1-\ve)^{p-2}t)-C.
$$
Hence, for all $t>0$
$$
at^{\gamma\beta}\ge (1-\ve)((1-\ve)^{p-2}t)^{\gamma\beta}\,b-C
$$
Letting $\ve, h\to 0$ and $C$ bounded we get $a\ge b$.
\qed

Here  is a  technical lemma that we have used.

\begin{lemma}\label{lemm.plcalc} If $f(x)=A|x|^{\gamma }$ with $0<\gamma< sp/(p-1)$, we have
$$
-\mathcal L_{s,p}f(x)= c(s,p,\gamma)A^{p-1}|x|^{\gamma(p-1)-sp}
$$
with $c>0$. For $\gamma\to sp/(p-1)$ then $-\mathcal L_{s,p}f(x)$ tends to a positive constant that can be infinite.
\end{lemma}

A sketch of the proof is as follows. The value of $-\mathcal L_{s,p}f(x)$ is finite for $x\ne 0$, the form of the function comes from scaling properties, the value of the constant comes from the limit as $x\to 0$. See in this respect the detailed study of Section 16.1 of \cite{VazFPL2-2020}.

\medskip

\noindent {\bf Eliminating the condition $sp<1$.} This restriction was only used at the beginning of the proof of Theorem \ref{thm.ss.growing} to ensure that there exists a finite minimal weak solution, a fact that depended on the existence on some kind of a priori bound from above, see Theorem \ref{thm.exist}. If $sp\ge 1$ we have to produce such an upper bound by a direct method. This is how we work.

1) We  modify  the initial function $f=|x|^\gamma$  near the origin so that $f_1\ge 1$, $f_1$ is increasing in $|x|$ and $C^2$ smooth, and the modification stops at $r=2$. Then we can assert that there is  a constant $C>0$ such that
$$
\mathcal L_{s,p}f_1(x)\ge -C f_1(x).
$$
Indeed,  on bounded sets this is clear, near infinity we can easily see that  $\mathcal L_{s,p}f_1(x)\sim \mathcal L_{s, p}f(x)$ so that $\mathcal L_{s,p}f_1(x)/ f_1\to 0$.

2) A candidate supersolution is then
$$
U(x,t)=K(t+1)^a f_1(x), \quad \mbox{with } a=1/(2-p).
$$
with $K$ large enough. Actually, since $a(p-1)=a-1$,
$$
U_t+\mathcal L_{s,p}(U)=aK(t+1)^{a-1}f_1+ K^{p-1}(t+1)^{a(p-1)} Lf_1
$$
$$\ge (t+T)^{a-1}f_1(x) (aK-CK^{p-1}).
$$
Therefore, if $K$ is large enough we have $aK-CK^{p-1}>0$ and we have proved that $U(x,t)$ is a  finite classical supersolution. Moreover,  $U(x,0)=K f_1(x)\ge f(x)=u_0(x)$. By the Maximum Principle applied to the approximations, we conclude that there is a finite solution in the limit and
$$
u(x,t)\le U(x,t).
$$
We thus have the a priori estimate for the solution to  power-like initial data with $0<\gamma<\gamma_1$. The rest of the proof of the theorem does not change. \qed

\subsection{Non existence of solutions for critical growth}\label{ssec.bu}

We now prove that the method of construction of solutions by approximation from below fails in the case of powers of the type $u_0(x)=C\,|x|^{\gamma}$ with $\gamma\ge sp/(p-1)$. We only need to examine the critical exponent
and we will find the phenomenon of \sl  instantaneous blow-up\rm, a typical occurrence in the evolution of linear or nonlinear parabolic equations with incompatible data. This phenomenon is well-known in the theory of semilinear of nonlinear heat equations with non-admissible data.

\begin{theorem}\label{thm.instbu} Let $u_0(x)=C\,|x|^{\gamma_1}$ with $\gamma_1= sp/(p-1)$, and let $u_{0n}$ be any nondecreasing sequence of nonnegative and integrable approximations to $u_0$. Let $u_n(x,t)$ be the corresponding semigroup solutions. Then for every $t>0$ and $x\in \ren$
\begin{equation}
u_n(x,t)\to \infty,
\end{equation}
and the convergence is  pointwise in $(x,t)$ and locally in $L^1_{loc}(\ren)$ for every $t>0$, and also in the local $L^1_{x,t}$ norm $Q_T$. Actually, we prove that
\begin{equation}
\lim_{\gamma\to\gamma_1}F_\gamma(x)\to \infty
\end{equation}
uniformly in $x$, where $F_\gamma$ is the self-similar solution constructed in Theorem \ref{thm.ss.growing} for exponents $\gamma<\gamma_1$.
\end{theorem}

\noindent {\sl Proof.~} (i) Let $\gamma<\gamma_1$ so that we now that the approximation process from below
converges in an increasing way towards a self-similar solution, $u_\gamma(x,t)=t^{\gamma\beta} \, F_\gamma(|x|\,t^{\beta})$ with $\beta(\gamma)>0$ as explained above.
Next, we observe that for every $\gamma<\gamma_1$ there is a constant $B_\gamma>0$ such that
$$
C\,|x|^{\gamma_1}\ge C\,|x|^{\gamma}-B_\gamma.
$$
Moreover, $B_\gamma\to 0 $ as $\gamma\to\gamma_1$. Using the invariance of the equation under vertical translations, we can find approximations $u_{0,n;\gamma}(x)$ to $C\,|x|^{\gamma}$ such that
$$
u_{0n}(x)\ge u_{0,n;\gamma}(x)-B_\gamma.
$$
so that $u_{n}(x,t)\ge u_{n;\gamma}(x,t)-B_\gamma.$ Passing now to the limit for $t=1$ we get
$$
\lim_{n\to\infty}u_n(x,1)\ge u_\gamma(x,1)-B_\gamma= F_\gamma(x)-B_\gamma
$$

(ii) Therefore, we are reduced to prove that the limit
$$
\lim_{\gamma\to\gamma_1}F_\gamma(x)=\infty \qquad \mbox{for every} \ x\ge 0.
$$
Since the $F_\gamma$ are radially symmetric and nondecreasing, so is the limit, hence we need to prove that
$$
\lim_{\gamma\to\gamma_1}F_\gamma(0)=\infty.
$$

(iii) Assume now that, at least along a subsequence ${\gamma_k\to\gamma_1}$ we have
$$
\lim_{\gamma\to\gamma_1}F_\gamma(1)=F_*(1)<\infty.
$$
Since we know that the function $G_\gamma(|x|)=F_\gamma(|x|)|x|^{-\gamma}$ is nonincreasing, the same happens when you pass to the limit $\gamma\to\gamma_1$, so that the limit function $F_*(x)$ is finite for every $|x|>1$ with
$$
F_*(x)\le F_*(1)|x|^{\gamma_1} \qquad \mbox{for } |x|>1,
$$
It is also finite and bounded for  $|x|\le 1 $ since $F_*(x) $ is monotone in $|x|$.

(iv) The monotonicity of $F_\gamma(|x|)|x|^{-\gamma}$ implies that there must be a limit
$$
\lim_{|x|\to\infty }F_*(x)|x|^{-\gamma}=A\ge 0.
$$
We will prove that $A\ge C$. We have for all $x$
$$
F_*(x)\ge F_\gamma(x)-B_\gamma\ge C\,|x|^{\gamma}-B_\gamma.
$$
Passing to the limit $\gamma\to\gamma_1$ we get $F_*(x)\ge C\,|x|^{\gamma},$ so that for $|x|\ge 1$ we have
\begin{equation}
C\,|x|^{\gamma_1}\le F_*(x)\le F_*(1)|x|^{\gamma_1}.
\end{equation}
The crucial observation is that since $\gamma_1(p-1)=sp$ this function obtains an infinite value when we apply  $\mathcal L_{s,p}$ to it.

(iv) With these estimates we pass to the limit in the weak formulation of the equation for the profile and
we arrive at a contradiction if $F_*(1)$ is finite. \nc That means that $F_*(1)=\infty$. The monotonicity of $G_{\gamma_1}$ implies then that $F_*(x)=\infty$ for every $x\ne 0$.

(v)  We still need to prove that $F_*(0)=\infty$.
We take the solution with same power data but constant $C_1>0$ smaller than $C$. Then we displace the origin
a unit distance to some $x_0$. we compare both initial data and we conclude that there is a $B>0$ such that
$$
C\,|x|^{\gamma_1}\ge C_1\,|x-x_0|^{\gamma_1}-B
$$
After doing the corresponding approximations and passing the limit we get
$$
F_*(0;C)\ge F_*(x_0,C_1)-B=\infty-B=\infty.
$$
The proof is complete. \qed

\begin{corollary} Instantaneous blow-up in the above sense happens for all solutions with locally integrable initial data $u_0\ge 0$ such that
$$
u_0(x)\ge A|x|^{\gamma_1}- f(x),
$$
where $f(x)$ is bounded or integrable or both.
\end{corollary}

We leave the proof as an exercise.


\section{Self-similar solutions with decaying data}\label{sec.exss2}

In this section we consider the case of initial power functions with negative exponents. Here we do not need the condition $sp<1$ since comparisons are done in the framework of the semigroup solutions of reference \cite{VazFPL2-2020}. We know that  self-similar solutions of the first type are impossible for exponents $\gamma<-sp/(2-p)$. The condition  $sp<1$ is not needed in this section. \nc

\subsection{Standard self-similarity theory}\label{sec.exss2.1}

 We start by the existence of self-similar solutions in the whole space that turn out to be bounded.

\begin{theorem}\label{decay.int} Let $1<p<2$ and $0<s<1$, and consider the Cauchy problem for the FPLE with initial data
$$
u_0(x)=A\,|x|^{\gamma}, \quad A>0, \ \gamma<0.
$$
Let $\gamma_2= \max\{-N, -sp/(2-p)\}<0$.

(i) For every $\gamma\in (\gamma_2, 0)$ and every $A>0$ there exists a unique self-similar solution of the first kind for this problem, of the form
$$
u(x, t)= t^{\gamma\beta} F(xt^{-\beta})
$$
with $\beta=1/(sp+\gamma(2-p))>0$. $F$ is continuous, positive and bounded.
$$
\gamma\beta=\frac{\gamma}{sp+\gamma(2-p)}<0.
$$

(ii) $\gamma<0$ implies that $F(y)$ is monotone decreasing in $r=|y|$. In that case $\|u(t)\|_\infty$ decreases like
$$
u(x,t)\le u(0,t)=F(0)\,t^{-\alpha}, \quad \alpha=|\gamma|\beta>\frac{|\gamma|}{sp}.
$$
The time exponent $\alpha=|\gamma|\beta>0$ increases with $|\gamma|$, as described below.

(iii) In all cases the profile is continuous and locally bounded. As $|x|\to \infty$ the profile decays in the precise way, $F(y)\sim A|y|^{\gamma}$.
\end{theorem}

\medskip
\begin{figure}[t!]
 \centering
 \vspace{-2cm}
 \includegraphics[scale=0.65]{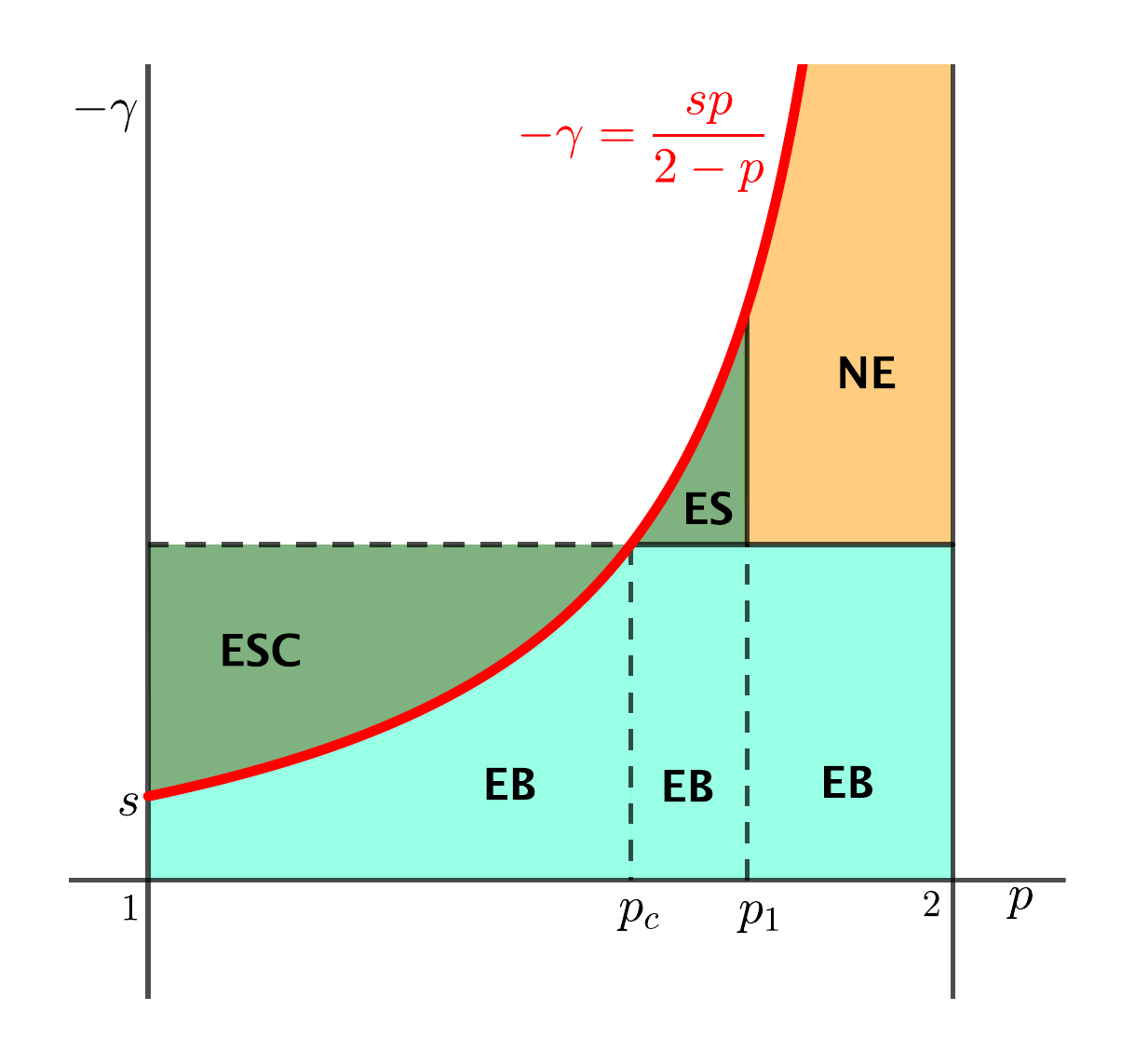}
 \caption{{Exponents of decaying self-similar solutions. Region EB represents existence of bounded solutions, ES denotes existence of expansive singular solutions.  ESC denotes existence of singular self-similar solutions of contractive type. 
 }}
  \label{fig:3}
\end{figure}

\noindent {\sl Proof of the Theorem.~}  (i)  The existence and uniqueness of a minimal weak solution
 uses the condition $\gamma>-N$ that ensures local integrability. The existence is a consequence of Section \ref{sec.ex} for $sp<1$, and is easy for $sp\ge 1$ arguing as follows. Even if the initial datum does not belong to any $L^p$ space, hence the standard semigroup theory does not apply, it is true that the truncated function
$$
(u_0(x)-C)_+\in L^q(\ren)
$$
for all $q<N/|\gamma|$ so that approximation from below by semigroup solutions is easy.

(ii) Self-similarity of the form \eqref{ext.sst1} is a consequence of the analysis of Section \ref{sec.self} once, thanks to the requirement  $sp+\gamma(2-p)>0$, which is the second condition satisfied by $\gamma$.
Radial symmetry immediately follows from uniqueness. Monotonicity in $r=|x|$ comes from Aleksandrov symmetry principle. The weak formulation for $F$ also follows.

(iii) Note that $\gamma_2=-N$ for $p_c\le p<2$, while $   \gamma_2=-sp/(2-p)<N$ for $1<p <p_c$.  We will see later that this range is optimal. \nc When $1<p \le p_c$ the values of $\beta(\gamma)$ cover the range $(1/(sp) ,\infty)$, and those of $\alpha=\beta|\gamma|$ go from 0 to infinity. On the other hand, when $p_c<p<2$ we get
$$
\frac1{sp}<\beta<\frac{p_c}{2N(p-p_c)}, \quad 0<\alpha<\frac{p_c}{2(p-p_c)}.
$$

(iv)  Because of the Lemma in the Appendix \nc we know that $u_t(x,0)\le 0$ and then we will have $u(x,t)\le u(x,0)$ for small $t$ and maximum principle will imply $u_t\le 0$ for all $t>0$.    \nc This also implies that  $F(y)\le A\,|y|^{\gamma}$, hence $F$ is bounded unless maybe at the origin.

(v) From $u_t\le 0$, i.e., $\mathcal L_{s,p}F(y)\ge 0$, and using equation \eqref{stateq.type1}, we have $\gamma F(y)\le yF'(y)$ so that
the function $G(r)=F(r)r^{-\gamma}$ is nondecreasing, hence it has a limit at infinity that must be
$$
\lim_{r\to\infty}F(r)r^{-\gamma}=A_1\le A.
$$
 The fact that $A_1=A$ happens because it $u(x,t)$ takes the initial data according to the conclusions of Theorem \eqref{thm.exist}. That immediately means that $F(y)$ is positive and bounded for all $y\ne 0$.

(iv)  We need to check that $F(0)$ is bounded. We use a new argument. If $p<p_c$ we consider the modified function
$$
u_{0,k}(x)=\max\{C, |x|^\gamma-k, 0\}
$$
This is an integrable function for which the semigroup theory applies. Since $|\gamma|<sp/(2-p)$ we can compare it with the explicit function of $U(x,t)$ of  \eqref{decay.vss.u2}. Hence, for $k$ large enough the
have
$$
u_k(x,t)\le (T-t)^{1/(2-p)} { |x|^{-sp/(2-p)} }.
$$
This means that at the time $t=T$ the solution $u_k(x,t)$ must be zero.  Next step is to observe that our self-similar solution can be compared with $u_k$ by vertical displacement. Indeed,
 $$
 t^{\gamma\beta} F(xt^{-\beta})\le u_k(x,t) +k
 $$
 For $t= T$ we have
 $$
 F(0)\le k\,T^{|\gamma|\beta}.
 $$

(v)  We need to check that $F(0)$ is bounded when $p\ge p_c$. We go back to the modified function that belongs to $L^1(\ren)$ so the smoothing effect shows that $u_k(x,t)$ if bounded in $\ren$ for all $t>0$. The end of the argument is the same. The case $p=p_c$ needs the smoothing effect for $u_0\in L^q$ with $q>1$.
\qed

\medskip

\noindent {\bf Smoothing effect.} The boundedness of the self-similar profile of these decaying self-similar solutions is also a consequence of the smoothing effect proved in \cite{BonSal2020} and \cite{DingMR4204564},  that we rephrase for our purposes as follows:

\begin{theorem}\label{smooth.thm}
Let $u$ be a weak solution of \eqref{frplap.eq} with $p>1$  corresponding to the initial datum $u_0\in L^{q_0}(\ren)$, with $q\ge 1$ when $p>p_c$ and $q>q_*=N(2-p)/sp$ when $1<p\le p_c$.  Then
\begin{equation}  \label{infty.bound}
\|u(\cdot,t)\|_\infty \leq   {C}\,{t}^{-\alpha}  \|u_0\|_{q_0}^{\mu}
\end{equation}
where $C=C(q,s,p,N)>0$  and the exponents come from dimensional consideration
$$
\mu  =     spq \beta(q), \qquad \alpha  = N \beta(q), \qquad \beta(q) =(N(p-2)+spq)^{-1}.
$$
\end{theorem}

The exponents come from dimensional considerations. Note that $\theta(q)>0$ when $q>q_*$ if $p\le p_c$.

\medskip

\noindent {\bf Remarks.} 1) Taking limits in the exponents we predict  for the linear case $p=2$ the values $\gamma_2=-N$, $\beta= 1/(2s)$ and $\alpha= \gamma/(2s)$. This is indeed the case, as proved in \cite{BonSirVazMR3614666}. Moreover, as $p\to 1$ the $\gamma $ interval stretches to $(-s,0)$.  We get $\beta\to 1/(s+\gamma)$ and $\alpha=\gamma/(s+\gamma).$ \nc

2) The explicit solution for $\gamma=-sp/(2-p) $ have extinction in finite time. Note  that expanding self-similar solutions of the first type are impossible for negative exponents such that $|\gamma|>sp/(2-p)$ if $p<p_c$. \nc

\medskip

\subsection{Singular self-similar solutions}\label{sec.sing.sss}

We examine here the existence of solutions in the region of first self-similarity but without the integrability condition. To be precise, we take $u_0(x)=|x|^{\gamma}$ with
$$
N\le -\gamma<sp/(2-p).
$$
 This region is non-empty when $p>p_c$. See Figure 3.

Remember that for all $p_c<p<p_1$  there is a very singular solution
\begin{equation}\label{decay.vss.u1.72}
V(x,t)=\,C_\infty  \,t^{1/(2-p)} { |x|^{-sp/(2-p)} }
\end{equation}

This is the first part of the result, that turns out to be positive.

\begin{theorem}\label{ex.sing.sss} Let $p_c<p<p_1$ and consider the Cauchy problem for the FPLE with initial data
$$
u_0(x)=A\,|x|^{\gamma}, \quad A>0,
$$

(i) For every $\gamma\in (-sp/(2-p),-N)$ and every $A>0$ there exists a unique self-similar solution of the first kind for this problem, obtained as minimal solution, and it has the form
$$
u(x, t)= t^{\gamma\beta} F(xt^{-\beta})
$$
with $\beta=1/(sp+\gamma(2-p))>0$. $F$ radially symmetric, nonincreasing along the radius $r=|x|$, it is positive for  $\ne 0$,  and has an isolated singularity at $r=0$. The time exponent is
$$
\gamma\beta=\frac{\gamma}{sp+\gamma(2-p)}<0.
$$
Now $\alpha=|\gamma|\beta>0$ increases with $|\gamma|$ with from $\alpha(-N) $ to infinity.

(ii) The profile decays in the precise way as $y\to\infty$, $F(y)\sim A|y|^{\gamma}$.

(iii) We have the equation for $F$ valid for all $x\ne0$

%
\end{theorem}

\noindent {\sl Proof of the Theorem.~} (i)   The existence comes from monotone approximations from below
in two steps. First, we prove existence for  the truncated solution
$$
 u_k(x,0)=(u_0(x)-k)_+\in L^1(\ren)
$$
by the construction of approximation with bounded functions, and using the bound from above by the VSS with large $T$. We have to check that the limit is a solution away from zero. Then we get a solution for the actual initial data by bounding the truncations from above the solutions
$$
v_k(x,t)=u_k(x,t)+k.
$$
Uniqueness follows as  a minimal weak solution.  Self-similarity of the form \eqref{ext.sst1} is a consequence of the analysis of Section \ref{sec.self} once we since $sp+\gamma(2-p)>0$.
Radial symmetry immediately follows from uniqueness. Monotonicity in $r=|x|$ comes from Aleksandrov symmetry principle. The weak formulation for $F$ also follows.

(ii) Because of the Lemma in the Appendix we know that $u_t(0)\le 0$ and then we will have $u(x,t)\le u(x,0)$ for small $t$ and maximum principle will imply $u_t\le 0$ for all $t>0$.
This also implies that  $F(y)\le A\,|y|^{\gamma}$, hence $F$ is bounded unless maybe at the origin.

(iii) it follows from $u_t\le 0$, i.e., $\mathcal L_{s,p}F(y)\ge 0$. Using equation \eqref{stateq.type1}, we have $\gamma F(y)\le yF'(y)$ so that
the function $G(r)=F(r)r^{-\gamma}$ is nondecreasing, hence it has a limit at infinity that must be
$$
\lim_{r\to\infty}F(r)r^{-\gamma}=A_1\le A.
$$
 The fact that $A_1=A$ happens because it $u(x,t)$ takes the initial data according to the conclusions of Theorem \eqref{thm.exist}. That immediately means that $F(y)$ is positive and bounded for all $y\ne 0$.

%

(iv) We show  $F(r)\to\infty$ as $r\to 0$. This is because the mass must be infinite. It is true if $sp<1$.
\qed

\medskip

\noindent {\bf Remark.} The existence of a finite  upper bound for $p_c<p<p_1$ is what produces the existence of the concept of minimal solutions obtained by approximation that will be weak solutions for $x\ne 0$.
The condition for $u^{p-1}$ integrable holds when $p<p_1$.

This remark is important because there is instantaneous blow-up when it fails.\nc

\begin{theorem}\label{nonex.sing.sss} Let $p_1\le p<2$ and consider the Cauchy problem for the FPLE with initial data
$$
u_0(x)=A\,|x|^{\gamma}, \quad A>0,
$$
with $-sp/(2-p)\le \gamma\le -N$. The above construction of a minimal solution blows up everywhere for $t>0$.
\end{theorem}

\noindent {\sl Proof.~} Here is a sketch of the proof.
The idea is that for $p_1\le p<2$  there is instantaneous blow up, since the VSS as limit of fundamental solutions is infinite and we can squeeze large mollified Dirac deltas below the initial data. Therefore, there is blow-up and there are no solutions. \qed

\subsection{Singular self-similar solutions of contractive type}\label{sec.sing.csss}

We have not studied the existence of solutions of the first kind in the remaining range of parameters $ \gamma<-sp/(2-p)$ where the analysis of Section predicts Subsection \ref{ssec.scaling.prd} predicts a solution with formula
$$
u(x, t)= t^{\gamma\beta} F(xt^{-\beta}).
$$
with $\beta$ negative and $\alpha=\gamma\beta>0$. We will not enter into the detailed analysis of this case here, we mention some details  because they provide interesting examples for next section.

Let us consider the nicest case,  where $1<p<p_c$ and $N>| \gamma|>-sp/(2-p)$. Then the initial data $u_0(x)=|x|^{\gamma}$ is a locally integrable function that also belongs to a Morrey space (see definition below) and the solution can be constructed as a self-similar solution of the first kind with contractive character in space, and no extinction in finite time (see next section). The singularity at $x=0$ is kept for all time in the form $u(x,t_1)\sim |x|^{\gamma}$ as $x\to 0$.

Further study of this topic is not needed here, It will appear elsewhere.


\section{ Study of the very fast diffusion range}\label{sec.vfd}

We consider solutions of the FPLE in the range $1<p<p_c$ with $N\ge 2$. There is no problem with the existence a unique strong solution for every initial datum $u_0\in L^q(\ren)$ with $1\le q<\infty$, and the semigroups are contractive in all $L^q$ spaces, including  the sup norm. We recall that for $p>p_c$ the law of mass conservation holds so there can be no extinction for any finite mass solution with $M>0$, and by comparison for any nonnegative solution with nontrivial initial data. The theory of finite-mass solutions was amply covered in \cite{VazFPL2-2020}.

\subsection{Study of extinction in finite time for $p<p_c$}

One of the main features of this range of parameters, $1<p<p_c$, is the possible loss of mass and eventual disappearance of the solution in a finite time. This is exemplified by the VSS \eqref{decay.vss.u2}, given by
$$
U(x,t)=\,C_\infty  \,(T-t)^{1/(2-p)} { |x|^{-sp/(2-p)} },
$$
where $p$ is any exponent in that range. An easy comparison implies the first extinction result that we will mention.

\begin{proposition} For every initial function such that
\begin{equation}\label{ext.eqn1}
u_0(x)\le A|x|^{-sp/(2-p)},
\end{equation}
there is extinction in finite time. More precisely, $T\le (A/C_\infty)^{2-p}$.\rm
\end{proposition}

Note that the function $|x|^{-sp/(2-p)}$ does not belong to any $L^q(\ren)$ space, but is a representative element of the Morrey space $M^{q*}(\ren)$ where $q_*=N(2-p)/sp$.

\smallskip

\noindent {\bf Definition. } \ For $q>1$ we define the Morrey space $\widetilde M^{q}(\ren)$ as the set of locally integrable $f$ such that
$$
 \int_{B_R(x_0)} |f(x)|\,dx \le C\,|R|^{N(q-1)/q} \quad \forall x_0\in \ren, R>0\,,
$$
with norm $\||f\||_q$ given by infimum of the $C$ in the above expression.

\smallskip

Note that  exponent $q_*$ is larger than 1 just for $p<p_c$. It grows to the value $N/s$ as $p\to 1$. The exponent is related to scale invariance and will appear in the next results.

\medskip

\noindent $\bullet$ A different sufficient condition is given by a result of Bonforte and Salort \cite{BonSal2020}.

\begin{proposition} {\rm \cite{BonSal2020}}\label{ext.bs}
Let $u$ be a weak solution of \eqref{frplap.eq} with $1<p<p_c$ corresponding to the initial datum $u_0\in L^{q_*}(\ren)$. Then, there exists a time $T=T(u_0)$  such that $u(t) \equiv 0$ for all $t\geq T$. Moreover, there is a constant $K(p,s)>0$ such that
\begin{equation}\label{ext.eqn2}
T\leq K(p,s)\|u_0\|_{q_*}^{2-p}.
\end{equation}
\end{proposition}

The space $L^{q_*}(\ren)$ includes in particular all bounded functions with compact support.
 Note that function \eqref{ext.eqn1} is not in $L^{q_*}(\ren)$, and yet it vanishes in finite time, so the space is not optimal in a sharp sense. It is though the border case for Lebesgue spaces, as we will see below. It is sometimes   called the Lebesgue {\sl extinction space}.

\medskip

\noindent $\bullet$ Let us now address the question of lower bounds for $T$, from which we can derive conditions for non-occurrence of extinction. In the case where $sp<1$ we can use the weighted $L^1$ estimate  to get the result.

\begin{proposition} let $1<p<p_c$ and $sp<1$. Let $u$ be a solution of \eqref{frplap.eq} with $1<p<p_c$ corresponding to a locally integrable initial datum $u_0\ge 0$. If $u$ has finite extinction time $T>0$, then $u(t)$ belongs to the Morrey space $M^{q*}(\ren)$  for all $T\ge t\ge0$. More precisely, the Morrey norm is bounded by
\begin{equation}\label{ext.eqn3}
\||u_0\||_{q_*}\le C(s,p)T^{1/(2-p)}.
\end{equation}
\end{proposition}

\noindent {\sl Proof.} Assuming the $\int u_0\,\var_R\,dx$ is finite and that the solution extinguishes at time $T>0 $, using in Corollary \ref{cor.1}  we get for all $0<t<T$  the estimate
\begin{equation}\label{cond.phi.Rext}
\left(\int_{B_R} u(t)\var_R\,dx\right)^{2-p} \le K(\varphi) R^{N(2-p)-sp}(T-t).
\end{equation}
It follows that the integral of $u$ in a ball of radius $R$ cannot grow more than $CR^{N(2-p)-sp}$ which is the definition of the Morrey space $\widetilde M^{q_*}(\ren)$. More precisely,
$$
\||u(t)\||_{q_*}\le K(\varphi)^{1/(2-p)}(T-t)^{1/(2-p)}.
$$
This proves the result.     \qed

\begin{corollary} Under the last assumptions, if  $u_0$ does not belong to $M^{q*}(\ren)$ then $T=\infty$. The following conditions imply positivity for all  time :

(i) $\lim_{|x|\to \infty} u(x)|x|^{sp/(2-p)}=+\infty$.

(ii)  $\lim_{|x|\to 0} u(x)|x-x_0|^{sp/(2-p)}=+\infty$ for some $x_0\in\ren$.

(iii) In particular, if the initial data is a power function,
$$
u_0(x)=|x|^{-|\gamma|}
$$
and $0<|\gamma|<N$, then $u(x,t)$ has finite extinction time only for $|\gamma|=\gamma_*=sp/(2-p)$.

This last result means that there are solutions with infinite-time blow-up in every Morrey space $\widetilde M^{q}(\ren)$ with $q\ne q_*$.
\end{corollary}

\noindent {\sl Proof.} Points (i) and (ii) are easy. As for (iii), note that (i) above includes the self-similar solutions with initial power function $u_0(x)=|x|^{-|\gamma|}$ with  $0<|\gamma|< \gamma_2=sp/(2-0)$. We know that each one is an example of function
in $\widetilde M^{q}(\ren)$ with $q>q_*$ that is positive for all  time.

On the hand,  case (ii) includes the solutions with initial power function $u_0(x)=|x|^{-|\gamma|}$ \ with  \ $N>|\gamma|> \gamma_2=sp/(2-p)$. They are not self-similar and they represent solutions with data in
$\widetilde M^{q}(\ren)$ with $1< q<q_*$ that do not vanish in finite time. \qed

From the  previous results we get the inclusions
\begin{equation}
L^{q_*}(\ren)\subset \mathcal{FET}(\ren)\subset \widetilde M^{q_*}(\ren),
\end{equation}
where $\mathcal{FET}$ is the set of solutions of equation \eqref{frplap.eq} with locally integrable data undergoing extinction in finite time. The second is proved for $sp<1$.
We conclude that $\widetilde M^{q_*}(\ren)$ is the optimal Morrey space for extinction.

\medskip

\noindent {\bf Question.}  \nc  What happens for $|\gamma |\ge N$ in this range?  The minimal solutions either  blow up
instantaneously or they are finite for $x\ne 0$. In any case the value at $x=0$ is infinite. \nc

\noindent {\bf Question.}  \nc What is the behaviour of solutions with data  $u_0\in \widetilde M^{q_*}(\ren)$, $u_0\not\in L^{q_*}(\ren)$?  We have a solution in that set that still vanishes in finite time, the VSS. Are there examples of solutions with $T=\infty$?
\nc

\medskip

\noindent {\bf Comment.} Questions of extinction for the FPLE posed in a bounded domain have been studied by authors like \cite{AbdABPeral2018}. For more information about the FPLE  in bounded domains see for instance
\cite{Mazon2016, Vazquez2016} and references therein.

\subsection{The question of mass decay}
We consider here solutions with finite mass and wonder what happens with the mass evolution. We know that for $p>p_c$ mass is conserved. The opposite happens here.

\begin{theorem}\label{thm.massto0} Let $1<p<p_c$ and let $u$ be a semigroup solution of the FPLE such that $u(x,t)\ge 0$ and it has finite initial mass, $\int_{\ren} u(x,0)\,dx=M>0$. Then,
$$
M_u(t):=\int_{\ren} u(x,t)\,dx  \to 0 \quad \mbox{as \ } \ t\to\infty.
$$
\end{theorem}

\noindent {\sl Proof.~}  Given $u_0(x)=u(x,0)\in L^1(\ren)$ we approximate it in $L^1$ by a sequence of bounded and integrable functions $u_{0,n}$. According to Proposition \ref{ext.bs} we know that the corresponding solution has a finite extinction $T_n$ time that depends on the norm of $u_{0,n}$ in $L^{q_*}(\ren)$. By the contraction property
$$
\|u(T_n)-u_n(T_n)\|_1\le \ve_n\to 0.
$$
But since $u_n(\cdot, T_n) \equiv 0$ we have $M_u(T_n)=\|u(\cdot, T_n)\|\le \ve_n$. Finally, we recall the basic property that the mass is a nonincreasing function of time.\qed

We already know that many of these integrable solutions extinguish in finite time. This is not the case for all of them. Thus, in the case $sp<1$ we know that solutions with integrable data $u_0$ that do not belong to the Morrey space
$\widetilde M^{q_*}(\ren)$ must vanish in infinite time. We do not know how to obtain these sharp results when $sp\ge 1$.  We can apply simple comparison and prove that $u_0\in L^1_{loc}$ and there are constants $C_n>0$
such thar $$
u_0(x)+C_n \ge n\,|x|^{-\gamma_*} \quad \forall n\ge 1,
$$
then there is no extinction in finite time.

\medskip

\noindent {\bf Exercise.} (i) Prove that solutions with behaviour near $x=0$ of the form
$$
u_0(x)\sim |x|^{\gamma}, \qquad \gamma>sp/(2-p),
$$
do not vanish in finite time, even if $sp\ge 1$.

(ii) Prove that solutions with behaviour as $|x|\to\infty$ of the form
$$
u_0(x)\sim |x|^{\gamma}, \qquad \gamma<sp/(2-p),
$$
do not vanish in finite time, even if $sp\ge 1$.

\noindent {\sl Hint.} You may use a construction based on sub-solutions and scaling.

\section{Conservation of mass for FPLE. Critical case}\label{sec.cmass.fpl}

 We have already shown that conservation of mass holds for strong semigroup solutions of the fractional $p$-Laplacian Equation with initial data in $L^1(\ren)$ in the range $p>p_c=2N/(N+s)$. This was  proved in \cite{Vazquez2020} for $p>2$ and in \cite{VazFPL2-2020} for $p_c<p<2$, and it was well-known for $p=2$. It is precisely the combined range where finite-mass self-similar solutions exist. It is not an easy proof, specially in the fast diffusion case. On the other hand, for $1<p<p_c$ (and $N\ge 2$) we have just proved that there are plenty of semigroup solutions with finite extinction time,  and moreover all finite-mass initial data give rise to solutions that lose all their mass with passing time.  Thus, we are only left with the question of mass conservation for the critical case $p=p_c$. Note that $1<p_c<2$ for all $N\ge 1$. Conservation does hold but it turned out to need a quite nontrivial proof.

\begin{theorem}\label{thm.mc.fple} Let $p=p_c$.  Let $u(x,t)\ge 0$ be the semigroup solution of the Cauchy Problem for the FPLE with initial datum $u_0\in L^1(\ren)$, $u_0\ge 0$.
Under the further assumption that $sp_c<1$, we can prove that for every $t>0$ the mass is conserved
$$
\int_{\ren} u(x,t)\,dx=\int_{\ren} u_0(x)\,dx.
$$
\end{theorem}

\noindent {\sl Proof.}    The proof we give here needs the extra assumption $sp_c<1$ to ensure that suitable duality methods work. Such an assumption was introduced in Section \ref{sec.wl1est}. It is always true for $N=1$, or for $N\ge 2$ and $s\le 1/2$. Note that  $sp_c=N(2-p_c)$.
Here is the full detail of the proof.

(i) As in similar proofs, we can make a reduction on the class of data. We may always assume that $u_0\in L^1(\ren)\cap L^\infty(\ren)$ and $u_0$ is compactly supported, say in the ball of radius $R_0$. If our form of mass conservation is proved under these assumptions, then it follows for all data $u_0\in L^1(\ren)$ as a consequence of the semigroup contraction property.

(ii) Under those extra assumptions, we can derive the condition of ``uniform small mass near infinity'', an important tool in what follows. We first recall the almost monotonicity in time, $(2-p)tu_t\le u$, that is due to scaling, see \cite{VazFPL2-2020}. We also recall two known facts:  if the initial support is contained in a ball, the solution is monotone in space (along outward cones), as follows from the Aleksandrov argument that we have used above, like in Section \ref{sec.masscon}. Using this and the uniform integrability, we conclude that for every constant $c>0$ there is a such that for large $r=|x|\ge 2R_0$ we have $u(x,t)\le cMr^{-N}$ for all $t>0$. Moreover,  we also know that the solution is uniformly positive before a possible extinction time  (a fact that we will exclude with our proof),  see argument above. We conclude that for any two times $0<t_1<t_2<T$ before the possible extinction, the mass of $u(t)$ outside a large ball of radius $R=R_\ve$  is less than $\ve$ for all $t\in (t_1,t_2)$, and  $R_\ve$ depends also on $t_1$ and $t_2$ but nothing more.

(iii) The technical part of the proof starts as the similar proofs in \cite{Vazquez2020, VazFPL2-2020} but it needs  a much more difficult technical treatment. First, we  recall that  $sp_c<N$, something that will  happen in all dimensions.
We do a calculation for  the tested mass. Taking a smooth and compactly supported test function $\varphi(x)\ge0 $, we have for $t_2>t_1>0$:
\begin{equation}\label{mass.calc1}
\left\{\begin{array}{l}
\displaystyle \left|\int u(t_1)\varphi\,dx-u(t_2)\varphi\,dx\right|\le \iiint
\left|\frac{\Phi(u(y,t)-u(x,t))(\varphi(y)-\varphi(x)}{|x-y|^{N+sp}}\right|\,dydxdt\\[10pt]
\le \displaystyle \left(\iiint  |u(y,t)-u(x,t)|^p\,d\mu(x,y)dt\right)^{\frac{p-1}{p}}
\left(\iiint |\varphi(y)-\varphi(x)|^p\,d\mu(x,y)dt\right)^{\frac{1}{p}}\,,
\end{array}\right.
\end{equation}
where $d\mu=|x-y|^{-N-sp}dxdy$. Space integrals are over $\ren$  and time integrals over $[t_1, t_2]$. We use the sequence of smooth test functions $\varphi_n(x)=\varphi(x/n)$ where $\varphi(x)$ is a cutoff function which equals 1 for $|x|\le 2$ and zero for $|x|\ge 3$. We take $n\ge 2R$. Then, have to consider different regions for the efficient calculation of the multiple integrals. Note that in \eqref{mass.calc1} we estimate integrals in absolute value (by taking absolute value of the integrand). In the next paragraphs we will forget the time integrals momentarily for ease of writing.

\medskip

\noindent $\bullet$ We first deal with the outer region $A_R=\{(x,y): |x|,|y|\ge R\}$. We have to estimate
\begin{equation*}\label{mass.calc2}
 \displaystyle I(A_R):= \iint_{A_n} \frac{|\Phi(u(y,t)-u(x,t))|\,|\varphi_n(y)-\varphi_n(x)|}{|x-y|^{N+sp}}\,dydx\,dt.
 \end{equation*}
 \nc
  We can split the absolute integral in two and calculate
$$
 \displaystyle I(A_R;1):= \iint_{A_R} u(x,t)^{p-1} \frac{\,|\varphi_n(y)-\varphi_n(x)|}{|x-y|^{N+sp}}\,dydx\,dt
$$
$$
\le \int_{A_R}  \,u(x,t)^{p-1}\, \mathcal M(\var_n)(x)\,dx\le \big(\int_{A_R}  \,u(x,t)dx\big)^{p-1}
\big(\int_{A_R} \mathcal M(\var_n)(x)^{1/(2-p)}\,dx\big)^{2-p}
$$
Now, the first integral is less than $\ve$ uniformly in $t$ by what was said above, while the second is bounded by
 $$
 \int_{A_R} (\mathcal M\var_n)(x)^{1/(2-p)}\,dx\le \int_{\ren} n^{-sp/(2-p)}(\mathcal M\var)(x/n)\,dx=
 \int_{\ren} (\mathcal M\var)(y)\,dy\le C.
 $$
 where $\mathcal M$  is the operator introduced in Section 2 of \cite{VazFPL3-2021}, by the formula
\begin{equation}\label{frplap.opM}
({\mathcal M}\var)(x):= P.V.\int_{\ren}\frac{|\var(x,t)-\var(y,t)|}{|x-y|^{N+sp}}\,dy\,,
\end{equation}
and has properties similar to $\mathcal L_{s,p}$ when $sp<1$.
 Here is where the assumption $sp<1$ is used.

 We proceed symmetrically with the integral $I(A_R;2)$ that uses  $ u(y,t)^{p-1}$ instead of $ u(x,t)^{p-1}$. Putting both together, we conclude that $I(A_R)\le C \ve $ where $C$ does not depend on $\ve$ or $n$.

\noindent $\bullet$  The argument in the inner region $B_n=\{(x,y): |x|,|y|\le 2n\}$ is simple. Since
$\varphi_n(x)-\varphi_n(y)=0$, we see that the contribution to the integral \eqref{mass.calc1} is zero.
Recall that we are using $n\ge 2R$.

 \noindent $\bullet$ We still have to make the analysis in  other regions in order to cover the whole domain $x,y\in \ren$. An option is to consider the two cross regions
$C_n=\{(x,y): |x|\ge 2n ,|y|\le R\}$ and $D_n=\{(x,y): |x|\le R ,|y|\ge 2n\}$. Both are similar so we will  only work out the contribution  in $D_n$. The idea is that we have an extra estimate: $|x-y|>n$ that avoids the singularity in the weight of the integrand. Forgetting again the time integrals   for the moment, we have $\var(x)=1$, $\var(y)=0$, $u(y,t)\le cMr^{-N}$ so that
$$
I(D_n)\le C(p)\,(I_1(D_n)+ I_2(D_n)),
$$
where
 \begin{equation*}
 \begin{array}{c}
  \displaystyle I_1(D_n)\le\iint_{D_n} |u(x,t)|^{p-1}\,d\mu(x,y)\le \\
\displaystyle  \big(\int_{B_R} \,|u(x,t)|^{p-1} dx
  \big(\int_{|x-y|>n} |x-y|^{-N-sp} dy\big)   \big)\\
 \displaystyle  \le Cn^{-sp} \int_{B_R} |u(x,t)|^{p-1}\,dx\,.
 \end{array}
\end{equation*}
Since $0<p-1<1$, we have
$$
\int_{B_R} |u(x,t)|^{p-1}dx\le R^{N(2-p)}\left(\int_{B_R} |u(x,t)|\,dx\right)^{p-1}
\le R^{N(2-p)}\|u(x,t_1)\|^{p-1}.
$$
Therefore, $ I_1(D_n)\le CR^{N(2-p)} n^{-sp} $ which tends to zero as $n\to\infty$ with a power rate.
As for $I_2(D_n)$, using the known space decay of the solution as $|y|\to \infty$ we have:
$$
 \displaystyle I_2(D_n)\le\iint_{D_n} |u(y,t)|^{p-1}\,d\mu(x,y)\le
  cM R^N\int_{|y|\ge 2n} |y|^{-N-sp}y^{-N(p-1)}\,dy
  $$
 where we have integrated first the $x$ variable. Then $I_2(D_n)\le Cn^{-sp-N(p-1)}R^N\,,$
which tends to zero as $n\to\infty$.

Same analysis for $I(C_n)$. Note that these regions overlap but that poses no problem. This concludes the proof.

\medskip

(iii) Going back to the formula \eqref{mass.calc1}, we conclude in the limit $n\to\infty$ that the mass is conserved for the whole time interval
$$
\int_\ren u(x,t)\,dx= \mbox{constant} \qquad \mbox{ for all } \ t_1\le t\le t_2.
$$
Since $t_1<t_2$ are arbitrary we conclude that mass is conserved,
$$
\int_\ren u(x,t)\,dx= \int_\ren u_0(x)\,dx \qquad \mbox{ for all } \  t\ge  0.
$$
In particular, there will be no extinction in finite time.
\qed

\medskip

\noindent {\bf Conservation of integral for signed solutions.} The above proof implies the conservation in time of the integral $\int u(x,t)\,dx$ for signed solutions of the FPLE, \bc with the some easy comparison arguments.\nc

\medskip
%
%
%

\noindent {\bf Open problem.} The problem of proving mass conservation for the solutions of the fractional PLE with critical exponent, treated in this section, is still open  if the restriction $sp<1$ does not hold. One might wonder if the proof available for $p>p_c$ can be adapted by some kind of continuity.


\section{Brief analysis of the  fast PME}

The first step in the investigation of weighted inequalities for fractional diffusion was done in
paper  \cite{BonfVaz2014MR3122168} where we studied the case of fractional porous medium equation
$$
u_t+(-\Delta)^s (u^m)=0
$$
 in the range of exponents $0<s<1,$ $0<m<1$. For that equation full duality can be applied so that the
 basic inequality, Theorem 2.2 of  \cite{BonfVaz2014MR3122168}, applies for two ordered solutions $u\ge v$ of  the  FPME.
Then, for all $0\le s,t <\infty$ we have
\begin{equation}\label{HP.fde}
\begin{split}
& \left(\int_{\ren}\big(u(t,x)- v(t,x)\big)\varphi_R(x)\dx\right)^{1-m}\le\\
&\left(\int_{\ren}\big(u(s,x)- v(s,x)\big)\varphi_R(x)\dx\right)^{1-m}
+ \frac{C_1 \,|t-s|}{R^{2s-N(1-m)}}
\end{split}
\end{equation}
with $C_1>0$ that depends only on $\var,m,N$\,. It is obtained without any further restriction on the range of $p$. For a time we saw no way of applying any duality  argument to the FPLE
until partial duality was obtained  here in Section \ref{sec.wl1est}.

The basic weighted $L^1$-estimate \eqref{HP.fde} allowed us to construct an existence and uniqueness theory for minimal solutions under a growth condition. In terms of powers it reads, $u_0(x)=O(|x|^\gamma))$ with $\gamma<\gamma_*={2s/(1-m)}$. It is similar to  what this paper contains in Section \ref{sec.ex}.

However, the study of self-similarity is missing in \cite{BonfVaz2014MR3122168}, so that the sharpness of the mentioned growth condition for existence could not be tested,  as we have done here for the FPLE with the proof of instantaneous blow-up in the critical power case $u_0(x)=A|x|^{\gamma_*}$.  This step can be done now for FPME.

For a basic theory of the FPME we refer to \cite{dP11, dP12}. For a survey papers see \cite{VazTNLD2016, MR3289370}, for more on fractional fast diffusion \cite{VazVolMR3192423}. The construction of the VSS for the FPME was done in \cite{VazqEurop2014}. They provide examples of singular self-similar solutions.

\section{Appendix}

Here is the proof of a technical monotonicity lemma needed in Section \ref{sec.exss}. It is based on abstract semigroup theory. There could be simpler proofs.

\begin{lemma}\label{lemma.tech} Let $u$ be a semigroup solution of the FPLE such that $u(x,0)=f(x)>0$ and $Lf=g\le af$ with $a>0$. Then, if $L$ is a nonlinear $(p-1)$-homogeneous operator with maximum principle, the  function $u(x,t)e^{cat}$ is nondecreasing in time for some constant $c>0$ that does blow up when $a\to 0$.
\end{lemma}

\noindent {\sl Proof.~} (i) We want to use the Crandall-Liggett theorem on semigroups generated by $m$-accretive operators in Banach spaces, the original reference is \cite{CrandallLiggett71}. We need $Lf$ to belong to the same Banach space. We take as first step the equation for implicit discretization
$$
f(x)+ hLf=v_0
$$
so that we have $v_0\le (1+ah)f$. In other words, the resolvent operator $R_h=(I+hL)^{-1}$ satisfies
$$
R_h(v_0)=f\ge (1+ah)^{-1}v_0,
$$
as a pointwise inequality in $\ren$. But we have by the homogeneity oh the operator $L$
$$
R_{h'}(\lambda w)=\lambda\,R_{h}(w), \quad h'=\lambda^{p-2}h.
$$
Applying this to the first inequality with $w=v_0$, $\lambda_1=(1+ah)^{-1}$ and
$$
h_1= \lambda_1^{p-2}h=(1+ah)^{2-p}h
$$ we get
$$
R_{h_1} R_h(v_0)\ge R_{h_1} ((1+ah)^{-1}v_0)= (1+ah)^{-1}R_{h}(v_0)\ge (1+ah)^{-2}v_0.
$$
Now with $\lambda_2= (1+ah)^{-2}$, $h_2=\lambda_2^{p-2}h$:
\[
R_{h_2}  (R_{h_1}R_h(v_0))\ge R_{h_2}((1+ah)^{-2}v_0)= (1+ah)^{-2}R_h(v_0)\ge (1+ah)^{-3}v_0.
\]
By iteration we get $h_n= \lambda_n^{p-2}h=(1+ah)^{n(2-p)}$, $\lambda_n=(1+ah)^{-n}$, and
$$
R_{h_n} \cdots (R_{h_1}R_h(v_0))\ge  (1+ah)^{-n}v_0
$$
When $h$ is very small and $nh\sim T$ we get $h_n\sim e^{aT(2-p)}h$,
$$
T_1=\sum_0^n h_j=h\sum_0^n (1+ah)^{j(2-p)}\sim c_1 h\frac{(1+ah)^{(n+1)(2-p)}-1}{ah}\sim c_2 (e^{aT(2-p)}-1)
$$
and $(1+ah)^{-n}\sim e^{-aT}$
so that the discretized $T_h(v_0;T)$ solution $S_h(v_0)$ behaves like
$$
S_h(v_0)(T_1)\ge  e^{-aT}v_0.
$$

\noindent (ii) We now use the Crandall-Liggett theorem. Since $v_{0,h}\to f$ as $h\to 0$, the discrete function of Step (i) constructed by iteration converges when $h\to 0$ and the number of steps goes to infinity in such a way that $\sum_1^n(h)=T_1$, and indeed it happens that
$$
\lim_{h\to 0} R_{h_n} \cdots (R_{h_1}R_h(v_0))=u(T_1)=S_{T_1}(f),
$$
the semigroup solution of the FPLE with initial data $f$ at time $T_1$. Taking the limit in the proved estimate
we get
$$
u(x,T_1)\ge f(x)\,e^{-aT}
$$
with $T_1\sim c_2 (e^{aT(2-p)}-1)$. For small $aT$ this relation is linear up to a constant. \qed

\medskip

\noindent {\bf Application to Theorem \eqref{thm.ss.growing}.} We want to apply that result to settle the monotonicity in time of the self-similar solutions treated in that result. To justify that, we approximate the initial data $u_0(x)=|x|^{\gamma}$. We know that
$$
\mathcal L_{s,p}u_0(x)= -c(s,p,\gamma)A^{p-1}|x|^{\gamma(p-1)-sp}\le 0
$$
with a singularity at $x=0$. We first cut the function at height $x=K>0$ and smooth to get some $u_{0k}$
that has a slightly larger $\mathcal L_{s,p}u_{0n}$ in the region before smoothing and small and integrable later.
We finally cut the lower part at a small height $a>0$ and that makes the new function $f$ to have a bounded
$\mathcal L_{s,p}f$ near the origin and smaller later. Therefore,
$$
\mathcal L_{s,p}f_{a,k}\le \ve
$$
Now raise the function by a constant $c_\ve$ so that $g=f+c_\ve$ satisfies
$$
\mathcal L_{s,p}g(x)\le \ve g
$$
Using the Lemma on $g$, we get almost monotonicity for $g$ with an exponential that depends only on $\ve$. Passing to the limit $\ve\to 0$ we get actual monotonicity in time for the self-similar solution $u(x,t)$. \qed

\medskip

\noindent {\bf Application to the decaying solutions of Section \ref{sec.exss2}}. The argument is similar but now the signs are reversed. \qed

\section{ Comments and open problems}

\noindent $\bullet$  There are several open questions when $s>1/2 $ and $sp\ge 1$. Are there any weighted $L^1$ inequalities? Is there an existence theorem for a sharp class of initial data as the one in Section \ref{thm.exist}?

\noindent $\bullet$ Even for $sp<1$, we may try to find a existence condition that is optimal, or at least a finer rate that the  growth rate close to $|x|^{sp/(p-1)}$ that we give in Section \ref{sec.ex}.  The power is sharp but we do not know to what extent the log conditions in Corollary \ref{cor.limitcond} are necessary. We recall that there is no $\la$ in the limit case $p=2$, i.e., for linear fractional diffusion.

\noindent $\bullet$  The theory of self-similar solutions can be done for anisotropic data of power type plus angle dependence. If we consider the form
$$
u_0(x)=A(\theta)|x|^{\gamma},   \qquad \theta=x/|x|,
$$
then under simple conditions on $A(\theta)$ like being positive and bounded above and below away from zero, the theory of Sections \ref{sec.exss}, \ref{sec.exss2} can be repeated, thus obtaining forward self-similar solutions with similarity profiles $F(\theta,y)$ that have less properties than the isotropic ones, but they are bounded above and below by them. Relaxing the requirements on $A(\theta)$ may lead to surprising new results.

\noindent $\bullet$ The existence theory for solutions with growing data in so-called superlinear case, $p>2$, has not been considered here. The study deserves its own space.

\noindent $\bullet$  We have not dealt with the systematic study of self-similar solutions of the second kind, and we give no information about possible solutions of the third type (with exponential time factors).

\noindent $\bullet$ For alternative definitions of the fractional $p$-Laplace operator we refer to \cite{DelTesGCVaz2020}. This a new field related to the numerical implementation.

\

\noindent {\textbf{\large \sc Acknowledgments.}} Author partially funded by Project  PGC2018-098440-B-I00 (Spain). Partially performed as an Honorary Professor at Univ. Complutense de Madrid.


\newpage

\noindent {\sc Address:}

\noindent Juan Luis V\'azquez. Departamento de Matem\'{a}ticas, \\ Universidad
Aut\'{o}noma de Madrid,\\ Campus de Cantoblanco, 28049 Madrid, Spain.  \\
e-mail address:~\texttt{juanluis.vazquez@uam.es}

\

\noindent {\bf Keywords: } Solutions with growing data, self-similar solutions, Nonlinear parabolic equations, $p$-Laplacian operator, fractional operators, extinction.

\medskip

\noindent {\bf 2020 Mathematics Subject Classification.}
  	35R11,   	
    35K55,  	
    35C06.       

\end{document}